\documentclass[reqno,11pt,a4]{amsart}  
\usepackage{amsmath}
\usepackage{amssymb} 
\usepackage{mathtools}   
\usepackage{graphicx}
\usepackage{graphics}  
\usepackage{color}
\usepackage{verbatim} 
\usepackage[normalem]{ulem} 
\usepackage[mathscr]{eucal}
\usepackage{upgreek}
\usepackage{enumerate} 
\usepackage{cite}
\usepackage{amsmath}
\usepackage{amsfonts}
\usepackage{amssymb}
\usepackage{bbm}
\usepackage{cancel}
\usepackage{graphicx}
\usepackage{multicol}
\usepackage[utf8]{inputenc}
\usepackage{framed}
\usepackage{setspace}
\usepackage{amsthm}
\usepackage{hyperref}
\usepackage{caption}
\usepackage{subcaption}
\usepackage{bbm}
\usepackage{graphicx}
\usepackage[titletoc]{appendix}
\numberwithin{equation}{section}  
 
\usepackage[refpage,noprefix]{nomencl}
\usepackage{nomencl} 

\newtheorem{theorem}{Theorem}[section] 
\newtheorem{lemma}[theorem]{Lemma} 
\newtheorem{proposition}[theorem] {Proposition} 
\newtheorem{cor}[theorem]  {Corollary} 
\newtheorem{remark}[theorem]  {Remark} 
\newtheorem{definition}[theorem] {Definition} 
 
\newtheorem{assump}[theorem]{Assumption}

\theoremstyle{definition}

 \usepackage{float}
\restylefloat{figure}

\DeclareMathAlphabet{\mathpzc}{OT1}{pzc}{m}{it}

\newcommand{\babs}[1]{{\bigl\lvert #1\bigr\rvert}}
\newcommand{\Babs}[1]{{\Bigl\lvert #1\Bigr\rvert}}

\DeclarePairedDelimiter{\abs}{\lvert}{\rvert}

\renewcommand{\L} {\Lambda} %

\def\d{\delta} 
\newcommand{\e} {\varepsilon}

\newfam\Bbbfam 
\font\tenBbb=msbm10 
\font\sevenBbb=msbm7 
\font\fiveBbb=msbm5 
\textfont\Bbbfam=\tenBbb 
\scriptfont\Bbbfam=\sevenBbb 
\scriptscriptfont\Bbbfam=\fiveBbb 
 
\newcommand{\B}     {\mathbb{B}} 
  
\newcommand{\R}     {\mathbb{R}} 
\newcommand{\Z}     {\mathbb{Z}} 
\newcommand{\N}     {\mathbb{N}} 
\renewcommand{\P}   {\mathbb{P}} 
\newcommand{\D}     {\mathbb{D}} 
\newcommand{\E}     {\mathbb{E}} 
\newcommand{\Q}     {\mathbb{Q}} 
 
 \newcommand{\floor}[1]{\left\lfloor #1 \right\rfloor}

\def\1{{\mathchoice {1\mskip-4mu\mathrm l}      
{1\mskip-4mu\mathrm l} 
{1\mskip-4.5mu\mathrm l} {1\mskip-5mu\mathrm l}}} 
 
\def\comment#1{} 
\newtheoremstyle{thm}{2ex}{2ex}{\itshape\rmfamily}{} 
{\bfseries\rmfamily}{}{1.7ex}{} 
 
\newtheoremstyle{rem}{1.3ex}{1.3ex}{\rmfamily}{} 
{\itshape\rmfamily}{}{1.5ex}{} 
 
\newcommand{\bB} {\boldsymbol{B}}

\newcommand{\bd} {\boldsymbol{d}}

\newcommand{\Fcal}   {{\mathcal F }}

\newcommand{\Mcal}   {{\mathcal M }} 
 
\newcommand{\Ocal}   {{\mathcal O }}

\newcommand{\Tcal}   {{\mathcal T }}

\def\dir{\mbox{dir}}

\newcommand{\LC}{\mathscr{L}\!\mathscr{C}}

 \newcommand{\ex}{{\rm e}} 
 
\renewcommand{\d}{{\rm d}}

\newcommand{\supp}{{\operatorname {supp}}}

\newcommand{\dist}{{\operatorname {dist}}} 
\newcommand{\diam}{{\operatorname {diam}}}

\newcommand{\tr}{{\operatorname {Tr}}}

\newcommand{\Exp}{\mathscr{E}\kern-0.2mm{\operatorname{xp}}}
\newcommand{\Log}{\mathscr{L}\kern-0.2mm{\operatorname{og}}}

\renewcommand{\emptyset} {\varnothing}

\newcommand\NoBlackBoxes{\global\overfullrule0pt}
\NoBlackBoxes

\newcommand\mycom[2]{\genfrac{}{}{0pt}{}{#1}{#2}}
\newcommand{\ek}[1]{\left[#1\right]}
\newcommand{\rk}[1]{\left(#1\right)}

\renewcommand{\dir}{{\mathrm{dir}}}

\newcommand{\df}{\delta^{\mathrm{free}}}
\newcommand{\Zf}{\mathrm{Z}^\mathrm{free}}

\renewcommand{\H}{\mathsf{H}}
\newcommand{\gfrak}{\mathfrak{g}}

\newcommand{\hk}[1]{^{(#1)}}

\newcommand{\exc}{\mathrm{exc}}

\renewcommand{\D}{\Delta}
\renewcommand{\gamma}{\delta^{\mathrm{exc}}}

\newcommand{\gk}[1]{\left\{#1\right\}}
\newcommand{\hc}[2]{^{\tiny{(#1),\,#2}}}

\newcommand{\vol}{\mathrm{S}}

\begin{document}

\title[\hfill Gibbs measures for the repulsive Bose gas\hfill]
{Gibbs measures for the repulsive Bose gas}


\thanks{}

\subjclass[2010]{Primary: 60K35; Secondary: 82B21; 82B41}

\keywords{Gibbs measures, Bose gas, Feynman representation}  
\begin{abstract}
 We prove the existence of Gibbs measures for the Feynman representation of the Bose gas with non-negative interaction in the grand-canonical ensemble. Our results are valid for all negative chemical potentials as well as slightly positive chemical potentials. We consider both the Gibbs property of marked points as well as a Markov--Gibbs property of paths.
        
\section{Introduction}
\subsection{The model}
In this paper, we prove the existence of Gibbs states for an ensemble of interacting Brownian loops in $\R^d$, with $d\ge 3$. The ensemble studied is also known as the \textit{Feynman representation} of the Bose gas, see Section \ref{subsectionFeyman} for background. A Brownian loop $\omega\colon[0,\beta j]\to\R^d$ is a continuous path with $\omega(0)=\omega(\beta j)$. In this work, the \textit{inverse temperature} $\beta>0$ is positive and $j$ is always a positive integer. We write $\ell(\omega)=j$ if $\omega$ is a loop of duration $\beta j$. 

The interaction between different loops is as follows: fix a weight function $\Phi\colon [0,\infty)\to [0,\infty]$. For two loops $\omega,\tilde \omega$, we set the \textit{pair interaction} $T$
\begin{equation}\label{EquationDefT}
    T(\omega,\tilde{\omega})=\sum_{n=0}^{\ell(\omega)-1}\sum_{m=0}^{\ell(\tilde\omega)-1}\int_0^\beta\Phi\rk{\abs{\omega(n\beta+s)-\tilde\omega(m\beta+s)}}\d s\, .
\end{equation}
The \textit{self-interaction} of a loop is given by
\begin{equation}\label{EquationDefW}
    W(\omega)=\frac{1}{2}\sum_{n=0}^{\ell(\omega)-1}\sum_{m=0}^{\ell(\omega)-1}\1\{n\neq m\}\int_0^\beta\Phi\rk{\abs{\omega(n\beta+s)-\omega(m\beta+s)}}\d s\, .
\end{equation}
For a collection of loops $I=\gk{\omega_1,\ldots,\omega_n}$, the total energy is then equal to
\begin{equation}
    \H(I)=\sum_{\omega\in I}W(\omega)+\frac{1}{2}\sum_{\mycom{\omega,\tilde\omega\in I}{\omega\neq\tilde\omega}}T(\omega,\tilde\omega)\, .
\end{equation}
Suppose we are given two collections of loops: $I_\L=\gk{\omega_1,\ldots,\omega_n}$ and $J_{\L^c}=\gk{\tilde{\omega_1},\tilde{\omega_2},\ldots}$. Here, assume that the loops in $I_\L$ are restricted (to be defined in the next section) to a bounded set $\L\subset\R^d$ and the loops in $J_{\L^c}$ are outside of $\L$. We can then define the \textit{Gibbs-kernel} $\delta_\L$ with boundary data $J_{\L^c}$ as
\begin{equation}\label{Eq:Intro:Kernel}
    \delta_\L(A|J_{\L^c})=\frac{1}{\mathrm{Z}_\L(J_{\L^c})}\int \1_A\{I_\L\cup J_{\L^c}\}\ex^{-\beta \H(I_\L)-\beta \sum_{\omega\in I_\L}\sum_{\tilde{\omega}\in J_{\L^c}}T(\omega,\tilde{\omega})+\beta\mu\sum_{\omega \in I_\L}\ell(\omega)}\d \rk{\omega_i}_{i=1}^n\, ,
\end{equation}
where $A$ is an event depending on the loops and $\mu\in\R$ is the \textit{chemical potential}. Here, $\mathrm{Z}_\L(J_{\L^c})$ is chosen such that $ \delta_\L(A|J_{\L^c})$ is a probability measure (in the first argument). The reference measure $\d \rk{\omega_i}_{i=1}^n$ will be specified in Equation \eqref{EquationML}.

Our main result is as follows: for a large class of weight functions $\Phi$, there exists a translation invariant measure $\gfrak$ on collections of interacting loops in $\R^d$, such that for every bounded function $F$ with compact local support, the Dobrushin–Lanford–Ruelle (DLR) equation holds:
\begin{equation}\label{IntroDLREquation}
    \int F(I_\L\cup J_{\L^c})\d\gfrak\rk{I_\L\cup J_{\L^c}}=\int \int F(I_\L\cup J_{\L^c})\delta_\L(\d I_\L|J_{\L^c})\d\gfrak\rk{J_{\L^c}}\quad\textnormal{for all }\L\subset\R^d\textnormal{ bounded}\, .
\end{equation}
This means that $\gfrak$ is a Gibbs measure with respect to the kernel $\rk{\delta_\L}_\L$.\\
The above equation is often abbreviated as
\begin{equation}
    \gfrak\delta_\L=\gfrak\, .
\end{equation}
as as well as slightly positive chemical potentials. We consider both the Gibbs property of marked points as well as a Markov--Gibbs property of paths.
\end{abstract}
\author{Tianyi Bai \and Quirin  Vogel}
\address[Tianyi Bai]{NYU Shanghai, 1555 Century Ave, Pudong, Shanghai, China, 200122}

\email{tianyi.bai@nyu.edu}
\address[Quirin  Vogel]{Technical University Munich (TUM), Boltzmannstrasse 3, Garching, Germany}
\email{quirin.vogel@tum.de}
 
\maketitle

\section{Introduction}
\subsection{The model}
In this paper, we prove the existence of Gibbs states for an ensemble of interacting Brownian loops in $\R^d$, with $d\ge 3$. The ensemble studied is also known as the \textit{Feynman representation} of the Bose gas, see Section \ref{subsectionFeyman} for background. A Brownian loop $\omega\colon[0,\beta j]\to\R^d$ is a continuous path with $\omega(0)=\omega(\beta j)$. In this work, the \textit{inverse temperature} $\beta>0$ is positive and $j$ is always a positive integer. We write $\ell(\omega)=j$ if $\omega$ is a loop of duration $\beta j$. 

The interaction between different loops is as follows: fix a weight function $\Phi\colon [0,\infty)\to [0,\infty]$. For two loops $\omega,\tilde \omega$, we set the \textit{pair interaction} $T$
\begin{equation}\label{EquationDefT}
    T(\omega,\tilde{\omega})=\sum_{n=0}^{\ell(\omega)-1}\sum_{m=0}^{\ell(\tilde\omega)-1}\int_0^\beta\Phi\rk{\abs{\omega(n\beta+s)-\tilde\omega(m\beta+s)}}\d s\, .
\end{equation}
The \textit{self-interaction} of a loop is given by
\begin{equation}\label{EquationDefW}
    W(\omega)=\frac{1}{2}\sum_{n=0}^{\ell(\omega)-1}\sum_{m=0}^{\ell(\omega)-1}\1\{n\neq m\}\int_0^\beta\Phi\rk{\abs{\omega(n\beta+s)-\omega(m\beta+s)}}\d s\, .
\end{equation}
For a collection of loops $I=\gk{\omega_1,\ldots,\omega_n}$, the total energy is then equal to
\begin{equation}
    \H(I)=\sum_{\omega\in I}W(\omega)+\frac{1}{2}\sum_{\mycom{\omega,\tilde\omega\in I}{\omega\neq\tilde\omega}}T(\omega,\tilde\omega)\, .
\end{equation}
Suppose we are given two collections of loops: $I_\L=\gk{\omega_1,\ldots,\omega_n}$ and $J_{\L^c}=\gk{\tilde{\omega_1},\tilde{\omega_2},\ldots}$. Here, assume that the loops in $I_\L$ are restricted (to be defined in the next section) to a bounded set $\L\subset\R^d$ and the loops in $J_{\L^c}$ are outside of $\L$. We can then define the \textit{Gibbs-kernel} $\delta_\L$ with boundary data $J_{\L^c}$ as
\begin{equation}\label{Eq:Intro:Kernel}
    \delta_\L(A|J_{\L^c})=\frac{1}{\mathrm{Z}_\L(J_{\L^c})}\int \1_A\{I_\L\cup J_{\L^c}\}\ex^{-\beta \H(I_\L)-\beta \sum_{\omega\in I_\L}\sum_{\tilde{\omega}\in J_{\L^c}}T(\omega,\tilde{\omega})+\beta\mu\sum_{\omega \in I_\L}\ell(\omega)}\d \rk{\omega_i}_{i=1}^n\, ,
\end{equation}
where $A$ is an event depending on the loops and $\mu\in\R$ is the \textit{chemical potential}. Here, $\mathrm{Z}_\L(J_{\L^c})$ is chosen such that $ \delta_\L(A|J_{\L^c})$ is a probability measure (in the first argument). The reference measure $\d \rk{\omega_i}_{i=1}^n$ will be specified in Equation \eqref{EquationML}.

Our main result is as follows: for a large class of weight functions $\Phi$, there exists a translation invariant measure $\gfrak$ on collections of interacting loops in $\R^d$, such that for every bounded function $F$ with compact local support, the Dobrushin–Lanford–Ruelle (DLR) equation holds:
\begin{equation}\label{IntroDLREquation}
    \int F(I_\L\cup J_{\L^c})\d\gfrak\rk{I_\L\cup J_{\L^c}}=\int \int F(I_\L\cup J_{\L^c})\delta_\L(\d I_\L|J_{\L^c})\d\gfrak\rk{J_{\L^c}}\quad\textnormal{for all }\L\subset\R^d\textnormal{ bounded}\, .
\end{equation}
This means that $\gfrak$ is a Gibbs measure with respect to the kernel $\rk{\delta_\L}_\L$.\\
The above equation is often abbreviated as
\begin{equation}
    \gfrak\delta_\L=\gfrak\, .
\end{equation}
\subsection{Gibbs property}
To make the concepts from the previous section more precise, we need to talk about \textit{local configurations}. Whereas in most statistical mechanics models, defining locality does not pose any problems, for our model this presents a big issue. The choice of locality is not purely cosmetic, as it dictates the definition of the Gibbs kernel $\rk{\delta_\L}_\L$. For each family of kernels, a distinct set of Gibbs measures may exist, see Section \ref{SubsectionOpenProblems}. Given a collection of loops $\{\omega_i\colon i\in I\}$ encoded in a point measure $\eta=\sum_{i\in I}\delta_{\omega_i}$, we give three ways to define the \textit{restriction} of $\eta$ to any set $\L\subset\R^d$:
\begin{itemize}
    \item The set $\eta_\L$ of loops \textit{started} inside $\L$. This point of view is most prominent in the mathematical literature, as it allows for the theory of decorated point processes to be applied. It corresponds to \textit{free} boundary conditions. See Figure \ref{fig:free} for an illustration.
\begin{figure}[h]
    \centering
    \includegraphics[width=0.5\linewidth]{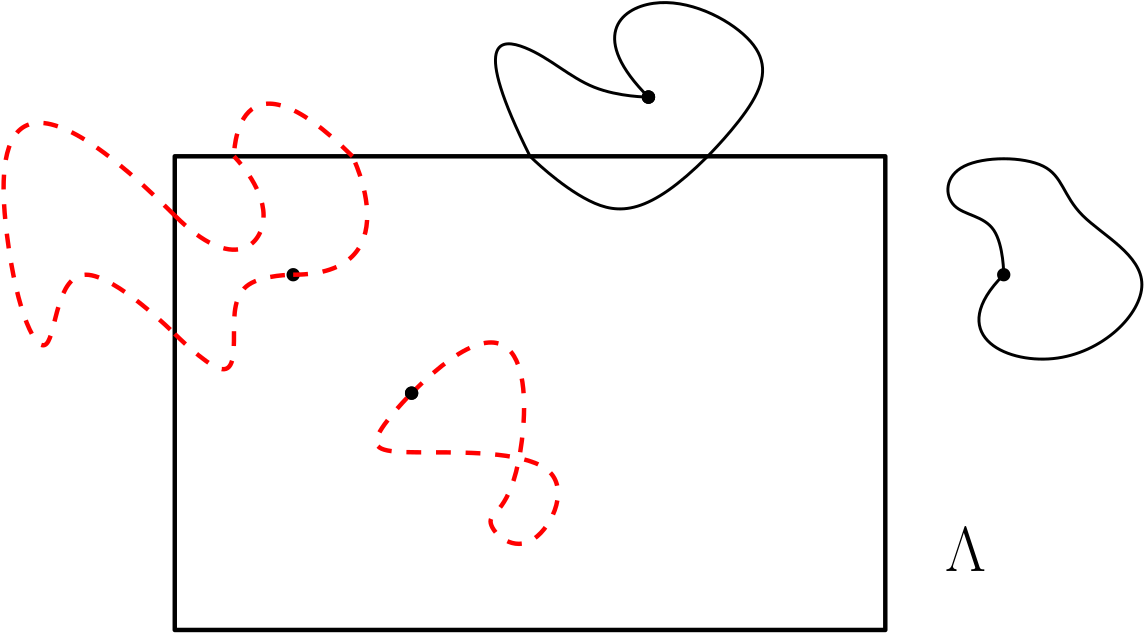}
    \caption{The set $\eta_\L$ in red (dashed), $\eta_\L^c$ in black. The Gibbs kernel $\df_\L$ (defined later) resamples the red loops.}
    \label{fig:free}
\end{figure}
    \item The set $\eta_\L^\dir$ of loops \textit{contained} in $\L$. This is the most natural definition in our setting, as it permits the definition of finite volume distributions for a wide range of choices $\Phi$. It corresponds to \textit{Dirichlet} boundary conditions. See Figure \ref{fig:dir} for an illustration.
\begin{figure}[h]
    \centering
    \includegraphics[width=0.5\linewidth]{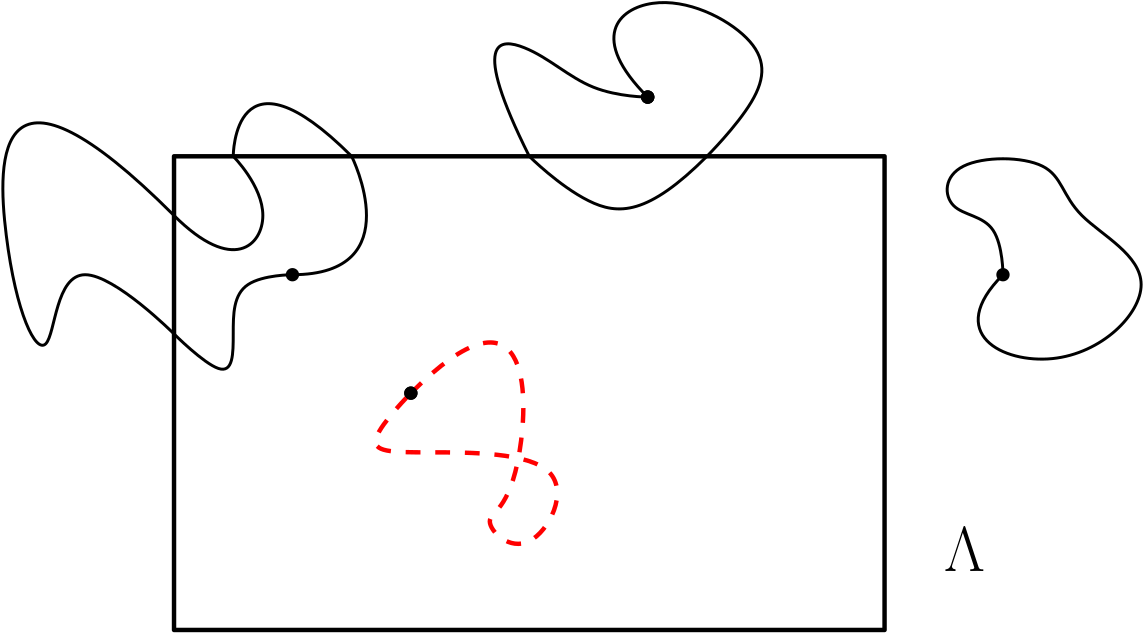}
    \caption{The set $\eta_\L^\dir$ in red (dashed), $\eta_\L^{\dir,c}$ in black. The Gibbs kernel $\delta_\L^\dir$ (defined later) resamples the red loops.}
    \label{fig:dir}
\end{figure}
    \item The set $\eta^\exc_\L$ of all \textit{paths} contained in $\L$. It consists of loops contained in $\L$ \textit{as well as} the excursions inside $\L$ of those loops which visit $\L^c$. This point of view is supported by the recent works connecting the Bose gas to random interlacements (see \cite{armendariz2021gaussian,vogel2020emergence,dickson2021formation}), where the other two notions of locality are no longer applicable. See Figure \ref{fig:exc} for an illustration.
\begin{figure}[h]
    \centering
    \includegraphics[width=0.5\linewidth]{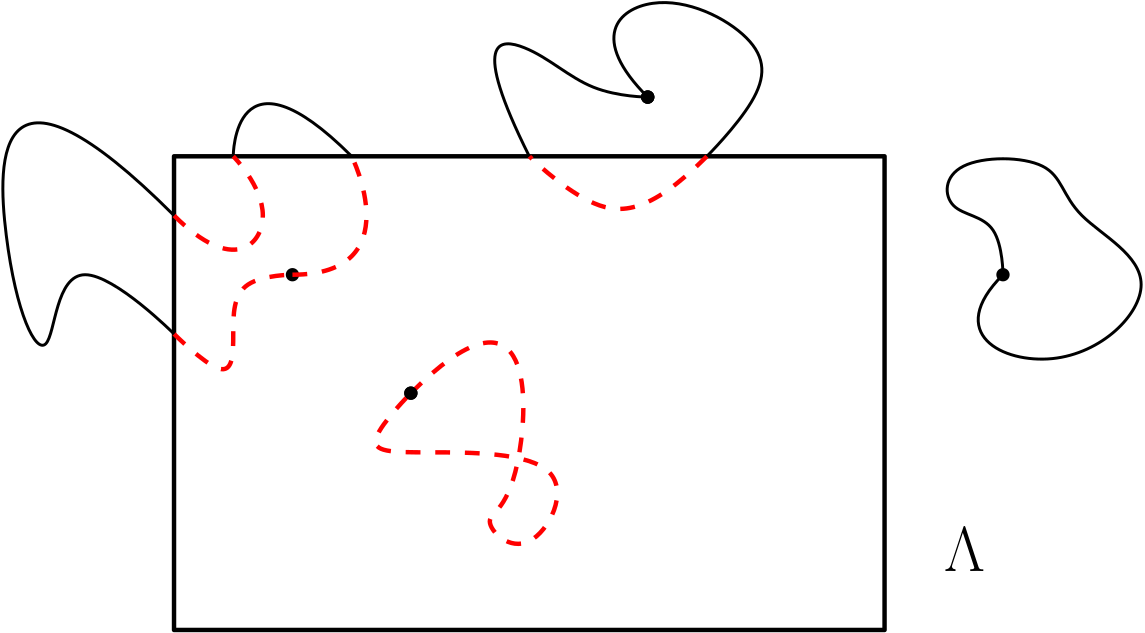}
    \caption{The set $\eta_\L^\exc$ in red (dashed), $\eta_{\L^c}^\exc$ in black. The Gibbs kernel $\gamma_\L$ (defined later) resamples the red paths.}
    \label{fig:exc}
\end{figure}
\end{itemize}
In our work, we consider the above three different families of kernels:
\begin{itemize}
    \item $\rk{\df_\L}_\L$, resampling the loops in $\eta_\L$.
    \item $\rk{\delta_\L^\dir}_\L$, resampling the loops in $\eta_\L^\dir$.
    \item $\rk{\gamma_\L}_\L$, resampling the paths in $\eta_\L^\exc$.
\end{itemize}
All the above kernels weigh configurations according to the weight $\ex^{-\beta \H}$, see Section \ref{SectionResults} for a rigorous definition.

The main result of our paper can be now made more precise: for $\Phi\ge 0$ satisfying a certain decay estimate,
\begin{center}
    there exists a probability measure $\gfrak$ which is Gibbs for all of the three kernels above.
\end{center}


%
%
\section{Results}\label{SectionResults}
Before stating the main result, we need some conditions on the interaction $\Phi$.
\begin{assump}\label{refassump}
Throughout the paper, we assume that $d\ge 3$. For the measurable weight function $\Phi\colon[0,\infty)\to [0,\infty]$, we assume that there exists $R>0$ and $\Psi\colon[0,\infty)\to[0,\infty)$ positive and decreasing with 
\begin{equation}\label{EquationBoundOnPsi}
    \int_R^\infty \Psi(x)x^{d-1}\d x<\infty,\quad\textnormal{and}\quad\Phi(x)\le \Psi(x)\quad \textnormal{for all}\quad x>R\, .
\end{equation}
Moreover, unless stated otherwise, all the domains $\L\subset\mathbb R^d$ are assumed to be connected, closed and satisfying the Poincar\'e cone condition, i.e. if on every point $x\in\partial\L$, there exists a cone $C$ with vertex $x$ such that $C\cap\L\cap\bB_x(r)=\{x\}$ for some $r>0$ small enough. Here, $\bB_x(r)$ is the ball centred at $x$ with radius $r$.
\end{assump}
Here, we have chosen to use $\Psi$ for the bounds on $\Phi$ far away from the origin. In most of the literature a separate function governs the behaviour of $\Phi$ close to the origin, see \cite{georgii1994large,ruelle1999statistical}. However, for our results the bound in Equation \eqref{EquationBoundOnPsi} suffices.\\
The Poincar\'e cone condition excludes domains with too irregular boundaries, which prevents paths to intersect at these boundaries. It is a purely technical condition, see Lemma \ref{lem:boundary_condition} for more.

Recall the chemical potential $\mu\in\R$ from Equation \eqref{Eq:Intro:Kernel}. Our main result is as follows:
\begin{theorem}\label{THM-Main}
Fix $\beta>0$. Under the Assumption \ref{refassump}, there exist a constant $c_\Phi>0$ such that for all $\mu\in\R$ with $\beta\mu<c_\Phi$, there exists a translation invariant probability measure $\gfrak=\gfrak(\beta,\mu,\Phi)$ on loop configurations such that in the sense of Equation \eqref{IntroDLREquation}
    \begin{equation}
        \gfrak\delta_\L^\dir=\gfrak\quad\textnormal{and}\quad\gfrak\df_\L=\gfrak\quad\textnormal{and}\quad\gfrak\gamma_\L=\gfrak\, ,
    \end{equation}
    for every bounded domain $\L\subset\R^d$, i.e., $\gfrak$ is a Gibbs measure for the three kernels above. See also Theorem \ref{theoremGibbsProperty} for a more precise restatement of the above result.

    A bound on $c_\Phi$ can be found in Equation \eqref{EquationcPhi}.
\end{theorem}
\begin{remark}
 If $\Phi\equiv 0$ Lebesgue almost everywhere, Theorem \ref{THM-Main} is trivial by standard Poisson theory in this case. It also holds true for all $\mu\in\R$. Hence, in the rest of the article we require that $\Phi\not\equiv 0$ and thus
    \begin{equation}
        \int_0^\infty\Phi(x)\d x>0\, .
    \end{equation}
\end{remark}
We now give a brief definition of the probability-kernels $\delta_\L^\dir$ and $\df_\L$. 

For $\eta=\sum_{\omega}\delta_\omega$, recall $\eta_\L$ and $\eta_\L^\dir$ in Figure \ref{fig:free}-\ref{fig:dir} (see also Equation \eqref{EquationDefinitionLocalConfig}). Set for $\L\subset\R^d$ bounded, 
\begin{equation}
    \mathrm{N}_\L(\eta)=\sum_{\omega\in\eta_\L}\ell(\omega)\, ,
\end{equation}
and
\begin{equation}
    \H_\L(\eta)=\sum_{\omega\in\eta_\L}W(\omega)+\frac 1 2\sum_{\omega\in \eta_\L}\sum_{\substack{\omega'\in \eta_\L\\ \omega'\ne\omega}}T(\omega,\omega')+\sum_{\omega\in \eta_\L}\sum_{\omega'\in\eta_\L^c}T(\omega,\omega')\, .
\end{equation}
In physical terms, one can think of $\mathrm{N}_\L$ as the number of particles in $\eta_\L$ and $\H_\L(\eta)$ as the interaction energy. 

Define $\P_{\L}$ the Poisson point process (see \cite[Chapter 24]{klenke2013probability} for a general definition) with intensity measure given by the \textit{Bosonic loop measure} $\mathrm{M}_{\L}$
\begin{equation}\label{EquationML}
    \mathrm{M}_{\L}=\sum_{j\ge 1}\frac{1}{j}\int_\L\d x\,\,\, \P_{x,x}^{\beta j}\, ,
\end{equation}
where $\P_{x,x}^{\beta j}$ is the unnormalized Brownian bridge measure from $x$ to $x$ in time $\beta j$ (see Equation \eqref{EquationBRidgeMeasure} for a definition). For $\L$ bounded, one can think (at least formally) of $\P_{\L}$ as
\begin{equation}
    \P_{\L}=\ex^{- \mathrm{M}_{\L}[\1]}\sum_{n\ge 0}\frac{ \mathrm{M}_{\L}^{\otimes n}}{n!}\, ,
\end{equation}
where $\1$ is the constant function.
We set for an event $A$
\begin{equation}
    \df_\L(A|\eta)=\frac{1}{\Zf_\L\rk{\eta_\L^c}}\int\1_A\{\xi+\eta_\L^c\}\ex^{-\beta\H_\L(\eta_\L^c+\xi)+\beta\mu \mathrm{N}_\L\rk{\xi}}\d\P_\L(\xi)\, ,
\end{equation}
and 
\begin{equation}
 \Zf_\L\rk{\eta_\L^c}=\int\ex^{-\beta\H_\L(\eta_\L^c+\xi)+\beta\mu \mathrm{N}_\L\rk{\xi}}\d\P_\L(\xi)\, .
\end{equation}
As $\beta>0$ and $\mu<c_\Phi/\beta$ remain fixed throughout the article, we do not include them in the notation of $\mathrm{M}_\L$, $\df_\L$ and $\P_\L$.

The measure $\P_\L^\dir$ is obtained from $\P_\L$ by restricting to Brownian motions contained in $\L$: $\P_\L^\dir$ is the Poisson point process with intensity measure $M^\dir_\L$, where $M^\dir_\L$ is given by $M_\L$ restricted to paths contained in $\L$, see Definition \ref{definitionRefProc}. We then define
\begin{equation}
    \delta_\L^\dir(A|\eta)=\frac{1}{\mathrm{Z}^\dir_\L\rk{\eta_\L^{\dir,c}}}\int\1_A\rk{\eta_\L^{\dir,c}+\xi}\ex^{-\beta\H_\L\rk{\eta_\L^{\dir,c}+\xi}+\beta\mu \mathrm{N}_\L\rk{\eta_\L^{\dir,c}+\xi}}\d\P_{\L,\beta}^\dir(\xi)\, ,
\end{equation}
where $\mathrm{Z}_\L^\dir\rk{\eta_\L^{\dir,c}}$ is the appropriate normalising constant such that $\delta_\L^\dir(A|\eta)$ is a probability measure in the first argument.

The definition of $\gamma_\L$ is significantly more involved, and we delay its definition to Equation \eqref{Equation deltaexc}.
Next, we give some properties of Gibbs measures with respect to our kernels.
\begin{proposition}\label{THMProp}
Suppose that Assumption \ref{refassump} holds. Let $\beta>0$, $\beta \mu<c_\Phi$.
\begin{enumerate}
    \item If $\mathsf{G}$ is Gibbs with respect to $\rk{\df_\L}_\L$ and $\alpha<c_\Phi-\beta\mu$, then for all functions $\psi$ on loops such that for every loop $\omega$, $\psi(\omega)\le \alpha \ell(\omega)$, we have
    \begin{equation}
        \mathsf{G}\left[\ex^{\sum_{\omega\in\eta_\D}\psi(\omega)}\right]<\infty\, ,\textnormal{ for all compact }\D\, .
    \end{equation}
    In particular, $\mathsf{G}[\ex^{\e \mathrm{N}_\L}]<\infty$, for $\e>0$ small enough.
    \item If $\mathsf{G}$ is Gibbs with respect to $\rk{\df_\L}_\L$, then for all $\alpha>0$
    \begin{equation}
        \mathsf{G}\ek{\rk{\sup_{\omega\in\eta_\D}\diam(\omega)}^\alpha}<\infty\, ,\textnormal{ for all compact }\D\, ,
    \end{equation}
    where $\diam(\omega)$ is the diameter of $\omega$, formally defined in Equation \eqref{Eq123three}.
\end{enumerate}
\end{proposition}
Similar (albeit more restrictive) statements can be given for $\mathsf{G}$ Gibbs with respect to $\rk{\delta^\dir_\L}_\L$. We leave this to the reader.
\subsection*{Structure of the paper}
In Section \ref{SectionDiscussion} we briefly introduce the Feynman representation of the partition function. We then comment on related literature and discuss the novelties in our approach. Finally, we point to some work in progress and open questions.

Section \ref{SectionConstructionOfGibbsMeasure} contains the proof of the main result, which can be furthermore split in several parts:
\begin{itemize}
    \item In Section \ref{SubsectionNotation} we introduce the notation.
    \item Next, in Section \ref{sec:4.2} we quantify the effects of the interaction $\H$ on a \textit{single loop}.
    \item In Section \ref{sec:kernels} we construct the different kernels. We furthermore prove that they form a consistent family. Approximations $\gfrak_n$ to the Gibbs measure $\gfrak$ are introduced.
    \item Section \ref{sec:4.4} introduces the specific entropy function $I$ and proves a bound for $I(\gfrak_n)$. This allows us to conclude that $\rk{\gfrak_n}_n$ has an accumulation point $\gfrak$.
    \item Section \ref{sec:4.5} is crucial: we show that the aforementioned convergence happens in a very fine topology.
    \item In the succeeding Section \ref{sec:4.6}, we introduce \textit{tempered configurations}.
    \item Finally, in Section \ref{sec:4.7} we prove the DLR equations for $\gfrak$, and show that they coincide for the different kernels. This establishes the main result Theorem \ref{THM-Main}.
    \item In addition, we give a proof for Proposition \ref{THMProp} in Section \ref{sec:4.8}.
\end{itemize}
In the Appendix, we provide a table with the frequently used notation. 
\section{Background and discussion}\label{SectionDiscussion}
\subsection{The Feynman representation}\label{subsectionFeyman}
Feynman in \cite{feynman1953atomic} used the theory of path integrals to give a stochastic representation of the Bose gas. For the purpose of giving context, we restrict ourselves to the partition function and refer the reader to Ginibre's notes (see \cite{ginibre1971some}) or the book by Bratteli and Robinson (see \cite{bratteli2003operator}) for an in-depth discussion. Furthermore, we introduce as little technical terms as possible. The complete definitions can be found in Section \ref{SectionConstructionOfGibbsMeasure}.

The \textit{partition function} $\mathrm{Z}_\L(\beta,\mu)$ of the \textit{grand-canonical} Bose gas in $\L\subset\R^d$ at \textit{inverse temperature} $\beta>0$ and \textit{chemical potential} $\mu\in\R$ is given by
\begin{equation}\label{EquationIntroFeynmanKac}
    \mathrm{Z}_\L(\beta,\mu)=\sum_{n\ge 0}\ex^{\beta \mu n}\tr_{\mathrm{L}^{2,+}(\L)^{\otimes n}}\rk{\ex^{-\beta\H_\L\hk{n}}}\, .
\end{equation}
Here, $\H_\L\hk{n}$ is the $n$-\textit{particle Hamiltonian} given by
\begin{equation}
    \H_\L\hk{n}=-\sum_{i=1}^n\Delta\hk{i}+\sum_{i,j=1}^n\Phi\rk{\abs{x_i-x_j}}\, ,
\end{equation}
where $\Delta\hk{i}$ is the standard Laplacian acting on the $i$-th coordinate and the second sum acts as a multiplication operator (here, $\Phi\colon[0,\infty)\to [0,\infty]$). The space $\mathrm{L}^{2,+}(\L)^{\otimes n}$ consists of those functions in $\mathrm{L}^{2}(\L)^{\otimes n}$ which are invariant under a permutation of coordinates, also called the \textit{Bosonic Fock space}. Feynman used the framework of what we call now the Wiener measure (rigorously constructed by Kac) to rewrite Equation \eqref{EquationIntroFeynmanKac} in terms of interacting trajectories. 
For the Feynman representation, we need a collection of loops $\omega_1,\ldots,\omega_m$, encoded in the point measure $\eta=\sum_{i=1}^m\delta_{\omega_i}$. 
Then Feynman's result reads
{\small\begin{equation}\label{EquationIntroPart}
     \mathrm{Z}_\L(\beta,\mu)=\sum_{n\ge 0}\sum_{j_1,\ldots,j_n}\!\!\!\frac{\ex^{\beta \mu (j_1+\ldots j_n)}}{n!\prod_{i=1}^n j_i}\!\int_{\L^n}\!\!\!\!\d x_1\ldots\d x_n\! \rk{\bigotimes_{i=1}^n \E_{x_i,x_i}^{\beta j_i}}\!\!\ek{\exp\gk{\!-\!\rk{\sum_{i=1}^nW(\omega_i)+\sum_{\mycom{j=1}{j\neq i}}^n T(\omega_i,\omega_j)}}} ,
\end{equation}}

\noindent where $W(\omega_i)$ gives the \textit{self-interaction} of each loop and $T(\omega_i,\omega_j)$ is the \textit{pair interaction} defined in the previous section. 

We can rewrite Equation \eqref{EquationIntroPart} with the help of the Poisson point process $\P_\L$:
\begin{equation}
     \mathrm{Z}_\L(\beta,\mu)=\ex^{\mathrm{M}_{\L}[1]}\E_{\L}\ek{\ex^{-\beta \H+\beta\mu \mathrm{N}_\L}}\, .
\end{equation}
Not only the partition function, but also particle density, correlation functions and other observables can be written in terms of the measure $\E_{\L}$, weighted by the factor of $\ex^{-\beta \H+\beta\mu \mathrm{N}_\L}$, see \cite{AV19}. This gives motivation to study the measure $\ex^{-\beta \H+\beta\mu \mathrm{N}_\L}\d\P_{\L} $.
\subsection{Literature}
The study of Gibbs measures has been pursued for many decades, going back to the works of Dobrushin, Lanford and Ruelle, see \cite{friedli2017statistical}. We only introduce a very small selection of references in the section, with a focus on those most relevant to our analysis. The most authoritative text on the subject is the monograph by Georgii \cite{georgii1988gibbs}, which considers lattice systems only. The study was complemented for particles positioned in $\R^d$ (see \cite{georgii1994large}) and for marked (or decorated) particles in \cite{georgii1993large}. A similar setting was used in the book on the subject by Preston, see \cite{preston2006random}. 

With regards to recent work, we highlight the very accessible paper by Dereudre \cite{dereudre2009existence}, in which the author considered geometry-dependent interaction between points in the plane $\R^2$. Recently, there has been interest in studying the existence of Gibbs measures for point processes decorated with random diffusion, see \cite{roelly2020marked} and \cite{zass2021gibbs}.

Besides proving the existence of Gibbs states, our work has another motivation: in \cite{adams2011variational}, the authors pose a minimisation problem over translation invariant probability measures in the context of an LDP result. According to the general Gibbs theory (see \cite[Chapter 15]{georgii1988gibbs} for example), it is conjectured that the measure $\gfrak$ is a solution for the minimisation problem. As this is technically rather involved, we have decided to prove that in separate publication.
%
\subsection{Novelties}
There exist several novelties in our proof of the existence of Gibbs measures. 

One novelty concerns the different Gibbs specifications (or kernels) used in this paper. While the specification of loops $\rk{\df_\L}_\L$ is standard in the literature on marked point processes, the other two specifications are not. The difference between $\rk{\df_\L}_\L$ and $\rk{\delta^\dir_\L}_\L$ is not too big conceptually, as only a surface order fraction of the loops exits the domain. However, for $\rk{\gamma_\L}_\L$, we have to work much harder as we need to separate the loops into excursions and paste them back together, which introduces additional dependencies. As mentioned previously, the kernel $\rk{\gamma_\L}_\L$ is motivated by the connection between random interlacements and the Bose gas, as studied in \cite{armendariz2021gaussian, vogel2020emergence,dickson2021formation}. Indeed, if the interaction $\Phi=0$ is set to zero everywhere, it was observed in \cite{vogel2020emergence} that the resulting superposition of loops and interlacement is Gibbs with respect to the resampling of loops/excursions inside a domain. There also exists more than one Gibbs measure for the kernel, see the next section for more. In the interacting case, this requires more work.
Furthermore, we would like to point out that the kernel $\rk{\gamma_\L}_\L$ encapsulates the notion of \textit{locality} in ``the most canonical'' way, as it limits itself strictly to all paths inside the domain $\L$.

Another important novelty is the \textit{absence of exponential integrability}. When one considers marked point processes with reference process $\mathsf{P}$, one usually\footnote{See for example  \cite{georgii1993large,dereudre2009existence,roelly2020marked,zass2021gibbs}} assumes that for some \textit{relevant observable} $\psi$
\begin{equation}\label{EquationIntroExpInt}
    \mathsf{P}\ek{\ex^{a\rk{\sum_{i}\psi(X_i)}}}<\infty \quad\text{for all }a\in \R\, .
\end{equation}
Here, our generic point process is written as $\sum_{i}\delta_{(i,X_i)}$ where $i$ is the location at which the mark $X_i$ is found. For example, in the Poisson Boolean model (see \cite{duminil2018subcritical}), $i$ is a point in $\R^d$ and $X_i$ is a ball centred at $i$. The function $\psi$ which maps $X_i\mapsto [0,\infty)$ encodes a characteristic of $X_i$ most relevant to the analysis. In the Poisson Boolean model, $\psi$ is usually a function of the radius of $X_i$. In the analysis of the of the Bose gas, there are two relevant observables:
\begin{enumerate}
    \item The \textit{particle number} $\ell(\omega)$. A loop parameterised on $[0,\beta j]$ is said to have $j$ particles, with $j\in \N$.
    \item The \textit{diameter} $\diam(\omega)=\sup_{0\le s,t\le \beta\ell(\omega)}\abs{\omega(s)-\omega(t)}$. 
\end{enumerate}
Controlling the particle number and the largest diameter for a group of loops is one of the main challenges of the proof. Usually, the exponential integrability from Equation \eqref{EquationIntroExpInt} helps with that. However, in our case
\begin{equation}
    \mathrm{M}_{\L}\ek{\ex^{a \ell(\omega)}}=\infty\, ,
\end{equation}
for all $a>0$. The same holds for the diameter. In that way our model is different from the aforementioned references.

We circumvent the above problem in a two step approach: we introduce an \textit{intermediate} measure $\tilde\P_\L$, which, while still being Poissonian, is in its decay properties much closer to the measure weighted by $\ex^{-\beta\H}$, see Equation \eqref{EquationIntermediateMeause}. We also make extensive use of the FKG-inequality and stochastic domination, see Lemma \ref{lem:FKG}. Translating the behaviour of the observables $\ell, S$ under $\gfrak$ into statements about topologies, we are able to circumvent the exponential integrability condition. We strongly believe that this approach has merits beyond the Bose gas model. Note that the papers proving the existence of Gibbs measures for marked point process with random paths (see \cite{zass2021gibbs} for example) usually use (super-)exponential integrability conditions for the diameter.

We also mention that we can handle (slightly) positive \textit{chemical potential} $\mu$. This is a first step into the direction of considering non-negative potentials (which is still open) and we believe that it should be possible to extend our proof for more cases. Indeed, we can show the existence of an accumulation point for $\gfrak_n$ for all \textit{superstable, regular} potentials, see \cite{georgii1994large} and our Remark \ref{remarksuperstable}. However, as we are motivated by the variational problem posed by \cite{adams2011variational}, we restrict ourselves to their setting: the interaction $\Phi$ is non-negative. 
\subsection{Open problems}\label{SubsectionOpenProblems} Having settled the question of existence, one may wish to discuss \textit{uniqueness} of Gibbs measures. We predict that this depends heavily on the kernel. For $\rk{\df_\L}_\L$, the theory of \textit{disagreement percolation} (see \cite{hofer2019disagreement} for example) is applicable for $\mu\ll 0$, at least in the case of a positive interaction. However, for the $\rk{\gamma_\L}_\L$ kernel, the problem is harder: in \cite{vogel2020emergence} it was shown that for $\Phi\equiv 0$ and $\mu,\tilde\mu\le 0$, the two measures
\begin{equation}
    \P_{\R^d,\beta,\mu}\quad\text{and}\quad \P_{\R^d,\beta,\tilde\mu}\otimes\P_\rho^\iota\, ,
\end{equation}
are both Gibbs with respect to the kernel $\rk{\gamma_\L}_\L$. Here, $\P_{\R^d,\beta,\mu}$ is the Poisson point process with intensity measure $\int_{\R^d}\d x\sum_{j\ge 1}\frac{\ex^{\beta\mu j}}{j}\P_{x,x}^{\beta j}$ and $\P_\rho^\iota$ is the Poisson point process of Brownian random interlacements with density $\rho\ge 0$, see \cite{sznitman2013scaling}. Note that $ \P_{\R^d,\beta,\mu}$ has particle density of at most $\beta^{-d/2}\zeta(d/2)$, while the particle density of $\P_{\R^d,\beta,\tilde\mu}\otimes\P_\rho^\iota$ can take any value. For more on this, we refer the reader to \cite{vogel2020emergence}.

Another open question related to this is whether it is possible to construct a version of the interlacements. This needs different tools, as the entropy argument breaks down.

The variational principle (see \cite{adams2011variational}) associated to $\rk{\df_\L}_\L$ has $\gfrak$ as the canonical candidate for its minimizer. This is the subject of further investigations.

\section{Proof}\label{SectionConstructionOfGibbsMeasure}

\subsection{Loop configurations}\label{SubsectionNotation}
In this section we introduce the notation used throughout the paper.

The basic objects of our analysis are Brownian loops $\omega$ of length $\beta j$ where $\beta>0$ is the inverse temperature and $j$ is a positive integer. Let 
\begin{equation}
    \Gamma_j=\gk{\omega\colon [0,\beta j]\to\R^d\text{ with }\omega(0)=\omega(\beta j)\text{ and }\omega\text{ continuous}},\,\Gamma=\bigcup_{j\ge 1}\Gamma_j\, .
\end{equation}
For any $\omega\in\Gamma_j$, we write $\ell(\omega)=j$. Write $\omega\in\L$ if $\omega(t)\in\L$ for all $t\in [0,\beta\ell(\omega)]$ and $\omega\cap\L\neq\emptyset$ if there exists $t\in [0,\beta\ell(\omega)]$ such that $\omega(t)\in\L$. The diameter is defined as
\begin{equation}\label{Eq123three}
    \diam(\omega)=\sup_{0\le s,t\le \beta \ell(\omega)}\abs{\omega(s)-\omega(t)}\, .
\end{equation}
Motivated by mathematical physics, we say that such $\omega$ \textit{represents $j$ particles} whenever $\omega\in\Gamma_j$. We equip $\Gamma$ with the topology induced by the topology of continuous functions on each $\Gamma_j$. The $\sigma$-algebra on $\Gamma$ is then taken as the associated Borel $\sigma$-algebra.

In this article, we study random point measures on $\Gamma$. Define $\Omega$ the space of all such point measures
\begin{equation}
    \Omega=\gk{\eta\colon \eta=\sum_{\omega\in I}\delta_\omega\text{ for }I\subset\Gamma\text{ at most countable}}\, .
\end{equation}
Equip $\Omega$ with the sigma-algebra of point measures $\Fcal$, as defined in \cite[Definition 24.1]{klenke2013probability}. We write $\omega\in\eta$ if $\eta(\omega)>0$. We furthermore write $\tau_x\eta$ for the configuration which is obtained from $\eta$ by shifting every loop by $x$, for $x\in\R^d$.

Our reference measures on $\Omega$ are given by Poisson point processes. Set $\P_x$ the measure of a $d$-dimensional Brownian motion, started at $x\in\R^d$. Let $\B_{x,y}^t$ be the Brownian bridge measure from $x$ to $y$ in time $t>0$, with $x,y\in\R^d$. Set 
\begin{equation}
p_t(x,y)=(2\pi t)^{-\frac d 2}\ex^{-\frac{\abs{x-y}^2}{2t}}\, ,
\end{equation}
the standard heat kernel of the Brownian motion. We define the (unnormalised) \textit{bridge measure}
\begin{equation}\label{EquationBRidgeMeasure}
    \P_{x,y}^t=p_t(x,y){\B_{x,y}^t}\, ,
\end{equation}
with total mass $p_t(x,y)$. Expectation with respect to $\P_{x,y}^t$ is denoted by $\E_{x,y}^t$. 
We will also need to introduce boundary conditions to our kernel. Define for an event $A$
\begin{equation}
    \B_{x,y}\hc{\L}{t}(A)=\B_{x,y}^{t}\rk{A|\omega\cap\L^c=\emptyset}\quad\text{and}\quad p_t\hk{\L}(x,y)=\d\P_x\rk{\omega(t)=\d y\text{ and }\omega\cap\L^c=\emptyset}\, ,
\end{equation}
for some domain $\L\subset\R^d$, and set
\begin{equation}
    \P_{x,y}\hc{\L}{t}=p_t\hk{\L}(x,y){\B_{x,y}\hc{\L}{t}}\, .
\end{equation}
Having established our measures on paths, we now define weights on $\Gamma$ and $\Omega$.
\begin{definition}\label{definitionRefProc}
Given a domain $\L\subset\R^d$ (bounded or unbounded) and an inverse temperature $\beta>0$, define
\begin{equation}\label{eq:M_L}
    \mathrm{M}_\L=\mathrm{M}_{\L,\beta}=\int_\L\sum_{j\ge 1}\frac{1}{j}\P_{x,x}^{\beta j}\d x\, ,
\end{equation}
and set $\P_\L=\P_{\L,\beta}$ the Poisson point process with intensity measure $\mathrm{M}_{\L,\beta}$. In particular, we have $\P_{\R^d}$ is the Poisson point process which samples loop on $\R^d$.

For $\L\subset\R^d$ bounded, we define \begin{equation}
    \mathrm{M}_\L^\dir=\mathrm{M}_{\L,\beta}^\dir=\int_\L\sum_{j\ge 1}\frac{1}{j}\P_{x,x}\hc{\L}{\beta j}\d x\, ,
\end{equation}
and set $\P_\L^\dir=\P_{\L,\beta}^\dir$ the Poisson point process with intensity measure $\mathrm{M}_{\L,\beta}^\dir$. 
\end{definition}
Notice that $\P_\L$ will only produce loops started within $\L$, while $\P_\L^\dir$ will only produce loops contained within $\L$. We want to reflect this notationally in our configuration. For this, we introduce for $\eta\in\Omega$
\begin{equation}\label{EquationDefinitionLocalConfig}
    \eta_\L=\sum_{\omega\in\eta}\delta_\omega\1\gk{\omega(0)\in\L}\quad\textnormal{and}\quad \eta_\L^\dir=\sum_{\omega\in\eta}\delta_\omega\1\gk{\omega\subset\L}\, ,
\end{equation}
see also the figures featured in the introduction.

We also define
\begin{equation}
    \eta^c_\L=\eta-\eta_\L\quad\text{and}\quad \eta_\L^{\dir,c}=\eta-\eta_\L^\dir\, .
\end{equation}
Furthermore, we define the sigma-algebras $\Fcal_\L$ and $\Fcal_\L^\dir$ induced by the projections $\eta\mapsto\eta_\L$ and $\eta\mapsto\eta_\L^\dir$. Also, similarly set $\Tcal_\L$ and $\Tcal_\L^\dir$ the sigma-algebras of $\eta\mapsto\eta_\L^c$ and $\eta\mapsto\eta_\L^{\dir,c}$.

In addition, notice that by \cite[Theorem 4.2.19]{karatzas2012brownian}, we avoid irregularity at boundaries of $\L$ given the Poincar\'e cone condition in Assumption \ref{refassump}:
\begin{lemma}\label{lem:boundary_condition}
Let $\L$ be a domain satisfying the Poincar\'e cone condition. Almost surely under $\P_{\mathbb R^d}$, loops are not starting at $\partial\L$, tangent to $\partial \L$, or intersecting (including self-intersections) at $\partial\L$.
\end{lemma}

Crucial are the following observables
\begin{definition}
For $\L\subset\R^d$ and $\eta\in\Omega$, we define
\begin{equation}\label{definitinvol}
    \mathrm{N}_\L(\eta)=\sum_{\omega\in\eta_\L}\ell(\omega)\quad\text{and}\quad \vol_\L(\eta)=\sup_{\omega\in \eta_{\L}}\diam(\omega)\, .
\end{equation}
\end{definition}

Recall the interaction terms $W$ and $T$ defined in the introduction.
\begin{definition}
For a bounded domain $\L$, we set
\begin{equation}\label{Equation H}
    \H_\L(\eta)=\sum_{\omega\in\eta_\L}W(\omega)+\frac 1 2\sum_{\omega\in \eta_\L}\sum_{\substack{\omega'\in \eta_\L\\ \omega'\ne\omega}}T(\omega,\omega')+\sum_{\omega\in \eta_\L}\sum_{\omega'\in\eta_\L^c}T(\omega,\omega')\, .
\end{equation}
We also set
\begin{equation}\label{Equation U}
    U(\eta;\xi)=\sum_{\omega\in\eta}\sum_{\substack{\omega'\in\xi\\\omega'\ne\omega}}T(\omega,\omega')\, ,
\end{equation}
for the interaction between two configurations.
\end{definition}
Due to the additivity of the energy, we have that:
\begin{lemma}\label{lem:H_trick_exc}
For any domains $\D\subset\L$ and any $\alpha,\beta\in\Omega$ not containing loops starting, tangent, or intersecting at $\partial \D\cup\partial\L$, 
if $\alpha-\alpha_\D=\beta-\beta_\D$, then
\begin{equation}\label{eq:Hamilton_shift}
\H_\L(\alpha)-\H_\D(\alpha)
=\H_\L(\beta)-\H_\D(\beta).
\end{equation}
In particular, we have \eqref{eq:Hamilton_shift} almost surely under $\P_{\R^d}$, see Lemma \ref{lem:boundary_condition}.
\end{lemma}

\subsection{Single-loop estimates}
\label{sec:4.2}
In this section, we discuss estimates based on a single loop $\omega$. We begin with a simple lemma, comparing the Brownian motion to the Brownian loop.
\begin{lemma}\label{LemmaCompareBridgeFree}
Suppose $G\colon\Gamma_j\to [0,\infty)$ bounded, such there exists $\e\in (0,\beta j)$,
\begin{equation}\label{eq:condition_bridge_to_bm}
    G\rk{\omega[0,\e]\oplus\omega_1[\e,\beta j]}= G\rk{\omega[0,\e]\oplus\omega_2[\e,\beta j]}\, ,
\end{equation}
for all $\omega,\omega_1,\omega_2$, then there exists $C=C_\e$ such that
\begin{equation}
    \E_{o,o}^{\beta j}[G]\le C\E_{o}\ek{G}\, .
\end{equation}
Here, the symbol $\oplus$ refers to the concatenation of two paths.
\end{lemma}
\begin{proof}
Set $\Fcal_t$ the sigma algebra generated by the projections $r\mapsto\omega(r)$ for $r\in [0,t]$.
Note that for $\tilde{G}$ an $\Fcal_t$-measurable function, we have that
\begin{equation}
    \E_{o,o}^{\beta j}\ek{\tilde{G}}=\E_o\ek{\tilde{G}(\beta j-t)^{-d/2}\ex^{-\abs{\omega(t)}^2/2(\beta j-t)}}\, .
\end{equation}
However, due to the condition on $G$, we can assume that it is $\Fcal_{\beta j-\e}$ measurable. Hence, the above remains bounded by a constant. This concludes the proof.
\end{proof}
Next we estimate the contribution of a single loop to the Hamiltonian.
\begin{lemma}\label{LemmaSingleLoopDecay}
There exists a constant $c=c_\Phi>0$ such that
\begin{equation}
    \E_{o,o}^{\beta j}\ek{\ex^{-\beta \H(\omega)}}=\Ocal\rk{\ex^{-c_\Phi j}},\,j\rightarrow\infty. 
\end{equation}
In fact, we can give the bound
\begin{equation}\label{EquationcPhi}
   3 c_\Phi\le -\log{\E\ek{\ex^{-\beta^2\int_0^1\Phi\rk{\beta^{1/2}\abs{B_s-B_1-\tilde{B}_s}}\d s}}}\, ,
\end{equation}
where $B,\tilde{B}$ are two independent Brownian motions started at the origin.
\end{lemma}
\begin{proof}
We abbreviate for $i=0,\ldots,j-1$,
\begin{equation}
    \omega\hk{i}=\rk{\omega(i\beta +t)}_{t\in [0,\beta]}\, .
\end{equation}
Then, by the non-negativity of the interaction
\begin{equation}
    \H(\omega)\ge \sum_{i=0}^{\floor{j/2}-1}T\rk{\omega\hk{2i},\omega\hk{2i+1}}\, .
\end{equation}
Further, by the previous equation and Lemma \ref{LemmaCompareBridgeFree},
\begin{equation}\label{eq:H_use_bridge_to_walk}
\begin{aligned}
    \E_{o,o}^{\beta j}\ek{\ex^{-\beta \H(\omega)}}
    &\le\E_{o,o}^{\beta j}\ek{\ex^{-\beta \sum_{i=0}^{\floor{j/2}-2}T\rk{\omega\hk{2i},\omega\hk{2i+1}}}}\\
    &\le C\E_{o}\ek{\ex^{-\beta\sum_{i=0}^{\floor{j/2}-2}T\rk{\omega\hk{2i},\omega\hk{2i+1}}}}
    =C\E_o\ek{\ex^{-\beta T\rk{\omega\hk{0},\omega\hk{1}}}}^{\floor{j/2}-2}
    \, .
\end{aligned}
\end{equation}
It remains to show that 
\begin{equation}
    \E_o\ek{\ex^{-\beta T\rk{\omega\hk{0},\omega\hk{1}}}}<1\, .
\end{equation}
Indeed, as $\int_0^\infty\Phi(x)>0$, there exists a set of $A$ of positive Lebesgue measure and a constant $\delta>0$, such that $\Phi(x)>\delta$ for all $x\in A$. Due to the absolute continuity of the finite dimensional distributions of the Brownian motion with respect to the Lebesgue measure, we have that
\begin{equation}
    \P_o\rk{\int_0^\beta \1\{\abs{\omega(s)-\omega(s+\beta)}\in A\}\d s>\e}>0\, .
\end{equation}
This concludes the proof of the first part. 

The lower bound on $c_\Phi$ follows for the scaling relation combined with the Markov property of the Brownian motion. Indeed,
\begin{equation}
\begin{aligned}
    \int_0^\beta\Phi\rk{\abs{\omega(s)-\omega(s+\beta)}}\d s
    &= \beta \int_0^1\Phi\rk{\abs{\omega(\beta s)-\omega(\beta s+\beta)}}\d s\\
    &=\beta \int_0^1\Phi\rk{\beta^{1/2}\abs{\omega( s)-\omega( s+1)}}\d s\, ,
\end{aligned}
\end{equation}
where the last equality is true in distribution. Furthermore,
\begin{equation}
    \omega(s)-\omega(s+1)=\omega(s)-\omega(1)-\rk{\omega(s+1)-\omega(1)}\, ,
\end{equation}
where $\omega(s+1)-\omega(1)$ is distributed like a standard Brownian motion and is independent from $\omega(s)-\omega(1)$, by the strong Markov property. This justifies \eqref{EquationcPhi}.
\end{proof}

\subsection{Gibbs kernels}\label{sec:kernels}
In this section, we introduce the Gibbs kernels $\rk{\delta^\dir_\L}_\L$, $\rk{\df_\L}_\L$ and $\rk{\gamma_\L}_\L$. We also introduce the approximations $\gfrak_n$ to the Gibbs measure $\gfrak$, and derive the FKQ-inequality (see Lemma \ref{lem:FKG}).
\subsubsection{The finite volume kernel}
We begin by defining the Gibbs kernel $\delta_\L^\dir$ for bounded domains $\L\subset\R^d$.
\begin{definition}
Define $\delta_\L^\dir\colon \Fcal\colon\Omega\to\R$ by
\begin{equation}\label{Equation delta}
    \delta_\L^\dir(A|\eta)=\frac{1}{\mathrm{Z}_\L^\dir\rk{\eta_\L^{\dir,c}}}\int\1_A\rk{\eta_\L^{\dir,c}+\xi}\ex^{-\beta\H_\L\rk{\eta_\L^{\dir,c}+\xi}+\beta\mu \mathrm{N}_\L\rk{\eta_\L^{\dir,c}+\xi}}\d\P_{\L,\beta}^\dir(\xi)\, ,
\end{equation}
with
\begin{equation}\label{Equation Z}
    \mathrm{Z}_\L^\dir\rk{\eta_\L^{\dir,c}}=\int \ex^{-\beta\H_\L\rk{\eta_\L^{\dir,c}+\xi}+\beta\mu \mathrm{N}_\L\rk{\eta_\L^{\dir,c}+\xi}}\d\P_{\L,\beta}^\dir(\xi)\, .
\end{equation}
\end{definition}

\begin{lemma}\label{LemmaFiniteIsFine}
The measure $\delta_\L^\dir$ is well-defined when  $\mathrm{N}_\L(\eta^{\dir,c}_\L),\H_\L(\eta^{\dir,c}_\L)<\infty$ and $\beta \mu<c_\Phi$.
\end{lemma}
\begin{proof}
In order that $\delta_\L$ is well-defined, we need to show that 
$0<\mathrm{Z}_\L<\infty$.
We have 
\begin{equation}
    \mathrm{Z}_\L^\dir\rk{\eta_\L^{\dir,c}}\ge \ex^{-\beta\H_\L\rk{\eta_\L^{\dir,c}}+\beta\mu \mathrm{N}_\L\rk{\eta_\L^{\dir,c}}}\P_{\L,\beta}^\dir(\xi=\emptyset)\, ,
\end{equation}
and hence $ \mathrm{Z}_\L^\dir\rk{\eta_\L^{\dir,c}}>0$ immediately by the condition $\H_\L(\eta^{\dir,c}_\L)<\infty$. 

Moreover, to prove $\mathrm{Z}_\L^\dir\rk{\eta_\L^{\dir,c}}<\infty$, it is enough to show that
\begin{equation}
 \mathrm{Z}_\L^\dir\rk{\eta_\L^{\dir,c}}\le \int \ex^{\sum_{\omega\in\xi}
\rk{-\beta\H\rk{\omega}+\beta\mu N\rk{\omega}}}\d\P_{\L,\beta}^\dir(\xi)<\infty\, ,
\end{equation}
due to the positivity of the interaction. By Campbell's formula (see 
\cite[Proposition 2.7]{last2017lectures}), this further reduces to 
\begin{equation}\label{Integrabilityofhamiltonian}
\int_\L \d x \sum_{j\ge 1}\frac{\ex^{\beta \mu j}}{j}\E_{x,x}^{(\L),\beta j}\ek{\ex^{-\beta\H(\omega)}}<\infty\,.
\end{equation}

Indeed, by Lemma \ref{LemmaSingleLoopDecay},
\begin{equation}
\begin{aligned}
\int_\L \d x \sum_{j\ge 1}\frac{\ex^{\beta \mu j}}{j}\E_{x,x}^{(\L),\beta j}\ek{\ex^{-\beta\H(\omega)}}
&\le\int_\L \d x \sum_{j\ge 1}\frac{\ex^{\beta \mu j}}{j}\E_{x,x}^{\beta j}\ek{\ex^{-\beta\H(\omega)}}\\
&=\Ocal\rk{|\L|\sum_{j\ge 1}\frac{\ex^{(\beta \mu-c_\Phi)j}}{j}}<\infty\,.
\end{aligned}
\end{equation}
This concludes the proof.
\end{proof}
Assume that from now on we have fixed some $\beta>0$ and $\beta \mu<c_\Phi$. Note that $\delta_\L^\dir(A|\cdot)$ is measurable with respect to $\Tcal_\L^{\dir,c}$ and that $\delta_\L^\dir(\cdot|\eta)$ is a probability measure for all $\eta$ with $\mathrm{Z}_\L(\eta_\L^{\dir,c})$ finite.
\begin{lemma}\label{LemmaConsistent}
The family $\rk{\delta^\dir_\L}_\L$ is a \textnormal{consistent} family, i.e. $\delta_\L^\dir\delta_\D^\dir=\delta_\L^\dir$ for $\D\subset\L$.
\end{lemma}
\begin{proof}
For simplicity, set $\mu=0$. The additional $\mu$-term does not affect the calculations, as it is linear.

Fix $\D\subset\L$, and by definition we have
\begin{equation}\label{eq:422}
    \delta_\L^\dir(\delta_\D^\dir(A|\cdot)|\eta)=\frac{1}{\mathrm{Z}_\L^\dir(\eta_\L^{\dir,c})}\int \d \P^\dir_{\L}(\xi)\delta_\D^\dir\left(A|\eta_\L^{\dir,c}+\xi\right)\ex^{-\beta\H_\L(\eta_\L^{\dir,c}+\xi)}\, .
\end{equation}
Notice that 
\begin{equation}
(\eta_\L^{\dir,c}+\xi)_\D^{\dir,c}=\eta_\L^{\dir,c}+\xi_\D^{\dir,c},
\end{equation}
then
\begin{equation}\label{eq:423}
    \delta_\D^\dir\left(A|\eta_\L^{\dir,c}+\xi\right)=
    \frac{1}{\mathrm{Z}_\D^\dir(\eta_\L^{\dir,c}+\xi^{\dir,c}_{\D})}\int\d \P^\dir_{\D}(\zeta)\1_A(\eta_\L^{\dir,c}+\xi_{\D}^{\dir,c}+\zeta)\ex^{-\beta\H_\D(\eta_\L^{\dir,c}+\xi_{\D}^{\dir,c}+\zeta)}\, .
\end{equation}
Note that $\xi=\xi_\L^\dir$, $\zeta=\zeta_\D^\dir$ since they are sampled from $\P^\dir_{\L}$ respectively $\P^\dir_{\D}$. By Lemma \ref{lem:H_trick_exc},
\begin{equation}\label{EquationHtrick}
    \H_\L(\eta_\L^{\dir,c}+\xi)-\H_\D(\eta_\L^{\dir,c}+\xi)=\H_\L(\eta_\L^{\dir,c}+\xi_{\D}^{\dir,c}+\zeta)-\H_\D(\eta_\L^{\dir,c}+\xi_{\D}^{\dir,c}+\zeta)\, .
\end{equation}
Therefore, we can rewrite $\delta_\L^\dir(\delta^\dir_\D(A|\cdot)|\eta){\mathrm{Z}^\dir_\L(\eta_\L^{\dir,c})}$ as
\begin{equation}
\begin{split}
    &\int\d \P^\dir_{\L}(\xi)\ex^{-\beta\H_\L(\eta_\L^{\dir,c}+\xi)}\frac{1}{\mathrm{Z}^\dir_\D(\eta_\L^{\dir,c}+\xi_{\D}^{\dir,c})} \int\d \P^\dir_{\D}(\zeta)\1_A(\eta_\L^{\dir,c}+\xi_{\D}^{\dir,c}+\zeta)\ex^{-\beta\H_\D(\eta_\L^{\dir,c}+\xi_{\D}^{\dir,c}+\zeta)}\\
    =&\int\d \P^\dir_{\L}(\xi)\ex^{-\beta\H_\D(\eta_\L^{\dir,c}+\xi)}\frac{1}{\mathrm{Z}^\dir_\D(\eta_\L^{\dir,c}+\xi_{\D}^{\dir,c})} \int\d \P^\dir_{\D}(\zeta)\1_A(\eta_\L^{\dir,c}+\xi_{\D}^{\dir,c}+\zeta)\ex^{-\beta\H_\L(\eta_\L^{\dir,c}+\xi_{\D}^{\dir,c}+\zeta)}\, .
\end{split}
\end{equation}
We can sample the loops contained in $\L$ by independently sampling loops contained in $\D$ and loops contained in $\L$ but not contained in $\D$, i.e., 
\begin{equation}\label{EquationPoissonDecom}
\d\P^\dir_{\L}(\xi)=\d\tilde\P_{\L\setminus\D}(\xi_{\D}^{\dir,c})\otimes\d\P^\dir_\D(\xi_\D^\dir)\, ,
\end{equation}
where $\tilde\P_{\L\setminus\D}$ is the Poisson point process with intensity measure 
\begin{equation}
    \int_\L\d x\sum_{j\ge 1}\frac{1}{j}\P_{x,x}^{\beta j}\circ\1\{\omega\subset \L\text{ and }\omega\not\subseteq\D\}\, .
\end{equation}
This allows us to rewrite $\delta_\L^\dir(\delta_\D^\dir(A|\cdot)|\eta){\mathrm{Z}^\dir_\L(\eta_\L^{\dir,c})}$ as
{\small
\begin{equation}
    \begin{aligned}
        &\int\d\tilde\P_{\L\setminus\D}(\xi_{\D}^{\dir,c})\int\d\P_\D^\dir(\xi_\D^\dir)\frac{\int\d \P^\dir_{\D}(\zeta)\1_A(\eta_\L^{\dir,c}+\xi^{\dir,c}_{\D}+\zeta)\ex^{-\beta\H_\D(\eta_\L^{\dir,c}+\xi_{\D}^{\dir,c}+\xi^\dir_\D)}\ex^{-\beta\H_\L(\eta_\L^{\dir,c}+\xi_{\D}^{\dir,c}+\zeta)}}{\mathrm{Z}^\dir_\D(\eta_\L^{\dir,c}+\xi_{\D}^{\dir,c})}\vspace{4mm} \\
        =&\int\d\tilde\P_{\L\setminus\D}(\xi_{\D}^{\dir,c})\int\d \P^\dir_{\D}(\zeta)\frac{\1_A(\eta_\L^{\dir,c}+\xi_{\D}^{\dir,c}+\zeta)}{\mathrm{Z}^\dir_\D(\eta_\L^{\dir,c}+\xi_{\D}^{\dir,c})}\ex^{-\beta\H_\L(\eta_\L^{\dir,c}+\xi_{\D}^{\dir,c}+\zeta)}\int\d\P_\D^\dir(\xi_\D^\dir) \ex^{-\beta\H_\D(\eta_\L^{\dir,c}+\xi_{\D}^{\dir,c}+\xi^\dir_\D)}\\
        =&\int\d\tilde\P_{\L\setminus\D}(\xi_{\D}^{\dir,c})\int\d \P^\dir_{\D}(\zeta){\1_A(\eta_\L^{\dir,c}+\xi_{\D}^{\dir,c}+\zeta)}\ex^{-\beta\H_\L(\eta_\L^{\dir,c}+\xi_{\D}^{\dir,c}+\zeta)}=\delta_\L^\dir(A|\eta){\mathrm{Z}_\L^\dir(\eta_\L^{\dir,c})}\, ,
    \end{aligned}
\end{equation}}
where we used that 
\begin{equation}
    \int\d\P_\D^\dir(\xi_\D^\dir) \ex^{-\beta\H_\D(\eta_\L^{\dir,c}+\xi_{\D}^{\dir,c}+\xi^\dir_\D)}=\mathrm{Z}_\D^\dir(\eta_\L^{\dir,c}+\xi_{\D}^{\dir,c})\, .
\end{equation}
This concludes the proof.
\end{proof}

We now construct another probability measure which only factors in self-interaction.
\begin{definition}\label{def:P^H_L}
Define
\begin{equation}
    \d\P_\L^\H(\eta)=\frac{1}{\E^\dir_\L\ek{\ex^{-\sum_{\omega\in\eta}(\beta\H(\omega)-\beta\mu\ell(\omega))}}}\ex^{-\sum_{\omega\in\eta}(\beta\H(\omega)-\beta\mu\ell(\omega))}\d\P_\L^\dir(\eta),
\end{equation}
then $\P_\L^\H$ is a Poisson point process with intensity measure $\mathrm{M}_\L^\H$ given by $\mathrm{M}_\L^\H[A]=\mathrm{M}_\L^\dir\ek{\ex^{-\beta\H+\beta\mu\ell}\1_A}$.
\end{definition}
The following lemma compares $\P_\L^\H$ to $\delta^\dir_\L(\cdot|\eta)$.
\begin{lemma}\label{lem:FKG}
We say a function $F\colon\Omega\to\R$ is \textnormal{increasing}, if
    \begin{equation}
        \supp(\eta_1)\subset\supp(\eta_2)\quad\Rightarrow\quad F(\eta_1)\le F(\eta_2)\, .
    \end{equation}
For every increasing function $F$,
\begin{equation}
    \delta^\dir_\L[F|\eta]\le \E_\L^\H\ek{F}\, .
\end{equation}
\end{lemma}
\begin{proof}
By definition,
\begin{equation}
    \delta^\dir_\L[F|\eta]=\frac{1}{\int \ex^{-\beta U(\xi;\xi)/2-\beta U(\xi;\eta_\L^c)}\d\P_\L^\H(\xi)}\int  F\rk{\xi}\ex^{-\beta U(\xi;\xi)/2-\beta U(\xi;\eta_\L^c)}\d\P_\L^\H(\xi)\, .
\end{equation}
Since $F$ is increasing and $\ex^{-\beta U(\xi;\xi)/2-\beta U(\xi;\eta_\L^c)}$ is decreasing,
by the FKG inequality for Poisson processes (see \cite[Theorem 20.4]{last2017lectures}), we have that
\begin{equation}
\int  F\rk{\xi}\ex^{-\beta U(\xi;\xi)/2-\beta U(\xi;\eta_\L^c)}\d\P_\L^\H(\xi)
\le \int  F\rk{\xi}\d\P_\L^\H(\xi)\int  \ex^{-\beta U(\xi;\xi)/2-\beta U(\xi;\eta_\L^c)}\d\P_\L^\H(\xi)\, ,
\end{equation}
which completes the proof.
\end{proof}

Heuristically, our Gibbs measure $\gfrak$ is the limit when domain $\L$ extends to $\mathbb R^d$, which motivates the following definition:
\begin{definition}
Let $\L_n$ be the centred cube of side-length $n$. We define
\begin{equation}
    \tilde{\gfrak}_n=\delta_{\L_n}^\dir(\,\cdot\,|\emptyset)\, .
\end{equation}
We can also extend our kernel periodically and set
\begin{equation}
    \hat\gfrak_n=\bigotimes_{x\in\Z^d} \delta_{xn+\L_n}^\dir(\,\cdot\,|\emptyset)\, ,
\end{equation}
as well as
\begin{equation}\label{Equation gn}
    \gfrak_n=\frac{1}{n^d}\sum_{x\in\L_n\cap\Z^d}\hat\gfrak_n\circ\tau_x\, .
\end{equation}
\end{definition}
Note that $\gfrak_n$ is translation invariant under shifts $\tau_x$ with $x\in\Z^d$, which will be important later on.
\begin{remark}
The reason that we choose the kernel $\rk{\delta^\dir_\L}_\L$ (instead of $\rk{\df_\L}_\L$ or $\rk{\gamma_\L}_\L$) for our approximation of the Gibbs measures is that this works for a wide range of weight functions $\Phi$. For instance, these definitions also work if $\Phi$ is a \textit{superstable potential}, see Remark \ref{remarksuperstable}. 
\end{remark}
Next, we show that the partition function is almost surely finite with respect to ${\gfrak}_n$
\begin{lemma}\label{Lem:PartFiniteFinite}
For all $n\ge 1$ and all $\Delta\subset\L_n$, we have that
\begin{equation}
    \gfrak_n\rk{\mathrm Z_\D^\dir=0}=0\, .
\end{equation}
\end{lemma}
\begin{proof}
   Note that by construction $\tilde{\gfrak}_n\rk{\H_{\L_n}+\mu\mathrm{N}_{\L_n}=\infty}=0$. Hence, by Lemma \ref{LemmaFiniteIsFine}, we have that $\mathrm Z_\D^\dir>0$ almost surely. 
\end{proof}

\subsubsection{The free kernel} 
Whereas our previous kernel $\delta^\dir_\L$ only resamples loops \textit{entirely} inside $\L$, we also need to introduce the \textit{free} kernel, which resamples all loops \textit{started} inside the domain. It is adapted to the topology $\rk{\Fcal_\L}_\L$ and resamples a larger class of loops. 
\begin{definition}
Define
\begin{equation}\label{Equation deltafree}
    \df_\L(A|\eta)=\frac{1}{\Zf_\L\rk{\eta_\L^c}}\int\1_A\{\xi+\eta_\L^c\}\ex^{-\beta\H_\L(\eta_\L^c+\xi)+\beta\mu \mathrm{N}_\L\rk{\xi}}\d\P_\L(\xi)\, ,
\end{equation}
where
\begin{equation}\label{Equation Zfree}
    \Zf_\L\rk{\eta_\L^c}=\int\ex^{-\beta\H_\L(\eta_\L^c+\xi)+\beta\mu \mathrm{N}_\L\rk{\xi}}\d\P_\L(\xi)\, .
\end{equation}
\end{definition}
Consistency of $\rk{\df_\L}_\L$ follows from the same argument in Lemma \ref{LemmaConsistent}. 

Moreover, one can easily check that Lemma \ref{lem:FKG} remains valid if we replace $\delta_\L^\dir$ by $\df_\L$.

Note that $ \Zf_\L\rk{\eta_\L^c}\ge \P_\L\rk{\xi=0}>0$, irrespective of $\eta$.
\subsubsection{The Gibbs--Markov kernel} In this section, we construct the Gibbs--Markov kernel. It is only relevant to the proof of Theorem \ref{TheoremPathDLR}, and this section can be skipped on first reading.

For each single $\omega\in\Gamma_j$ such that $\omega\cap\L\neq\emptyset$ and $\omega\cap\L^c\neq\emptyset$, we can split it into two parts, \begin{equation}
\omega^{\mathrm{exc}}_{\L}(t)=\omega(t)\quad\text{whenever}\quad\omega(t)\in\L^o,\,\text{otherwise undefined},
\end{equation}
where $\L^o$ is the interior of $\L$,
and 
\begin{equation}
\omega^{\mathrm{exc}}_{\L^c}(t)=\omega(t)\quad\text{whenever}\quad\omega(t)\in\L^c,\,\text{otherwise undefined}.
\end{equation}
Since $\L^o$ is open and $\omega$ is continuous, we know that the domain of definition for $\omega^{\mathrm{exc}}_\L$ is a (relatively) open subset of $[0,\beta j]$, thus a union of countably many intervals. Denote these intervals by $((a_i,b_i))_i$ (the interval touching $0$ and $\beta j$ may be closed), and we may define $\bd_\L(\omega)$ as the collection of their end-points and duration, 
\begin{equation}
\bd_\L(\omega)=\{(\omega(a_i),\omega(b_i),(b_i-a_i))_i\}\subset \mathbb R^d\times\mathbb R^d\times\mathbb R_+.
\end{equation}
Lift this construction to configurations, we define 
\begin{equation}
    \eta_{\L}^{\mathrm{exc}}=\sum_{\omega\in\eta}\delta_{\omega_\L^{\mathrm{exc}}},\quad
    \eta_{\L^c}^{\mathrm{exc}}=\sum_{\omega\in\eta}\delta_{\omega_{\L^c}^{\mathrm{exc}}},\quad
    \bd_\L(\eta)=\bigcup_{\omega\in\eta}\bd_\L(\omega)\, .
\end{equation}

By Lemma \ref{lem:boundary_condition}, $\P_{\mathbb R^d}$-almost surely, there is a unique way to reconstruct $\eta$ from $\eta_{\L}^{\mathrm{exc}}$ and $\eta_{\L^c}^{\mathrm{exc}}$ by gluing the paths at their endpoints. We abbreviate this procedure by
\begin{equation}\label{Equationexceta}
\eta=\eta_{\L}^{\mathrm{exc}}+\eta_{\L^c}^{\mathrm{exc}}.
\end{equation}

For $x,y\in\partial\L$ and $t>0$, define 
\begin{equation}
    \P_{x,y}\hc{\L,e}{t}\text{ the law of an excursion in }\L\text{ of length }t,\text{ from }x\text{ to }y\, .
\end{equation}
Given an at most countable collection $E=\rk{\rk{x_i,y_i,t_i}}_i$, we can sample a configuration $\zeta=\sum_i\delta_{f_i}$ so that 
\begin{equation}
    f_i \text{ is distributed by } \P_{x_i,y_i}\hc{\L,e}{t_i}\, .
\end{equation}
Denote this law by
\begin{equation}
    \P_\L^\mathrm{exc}\rk{\cdot\big|E}.
\end{equation}
In particular, we abbreviate
\begin{equation}
\P_\L^\mathrm{exc}\rk{\cdot\big|\bd_\L(\eta)}
=\P_\L^\mathrm{exc}\rk{\cdot\big|\eta}
=\P_\L^\mathrm{exc}\rk{\cdot\big|\eta^{\mathrm{exc}}_{\L^c}}.
\end{equation}

We can then define the following free kernel which resamples the dotted lines in Figure \ref{fig:exc}
\begin{equation}\label{EquationDefinitionOfQ}
    \Q_\L(A|\eta)=\int\1_A\{\eta_{\L^c}^{\mathrm{exc}}+\xi+\zeta\}\d\P_\L^\dir(\xi)\otimes\d\P_\L^\mathrm{exc}\rk{ \zeta\big|\eta}\, .
\end{equation}
\begin{lemma}\label{lem:consistency_Q}
For any $\D\subset\L$, $\P_{\mathbb R^d}$-almost surely any $\eta$ and any bounded and measurable function $f:\Omega\rightarrow\mathbb R$, (see Figure \ref{fig:Q_consistency})
\begin{equation}\label{eq:Q_consistency}
\int \d\Q_\L(\xi|\eta)f(\xi+\eta^{\mathrm{exc}}_{\L^c})=
\int \d\Q_\L(\xi|\eta)\int \d\Q_\D(\zeta|\xi+\eta^{\mathrm{exc}}_{\L^c})
f(\zeta+\xi^{\mathrm{exc}}_{\D^c}+\eta^{\mathrm{exc}}_{\L^c}).
\end{equation}
\end{lemma}
\begin{proof}
\begin{figure}
\includegraphics[width=10cm]{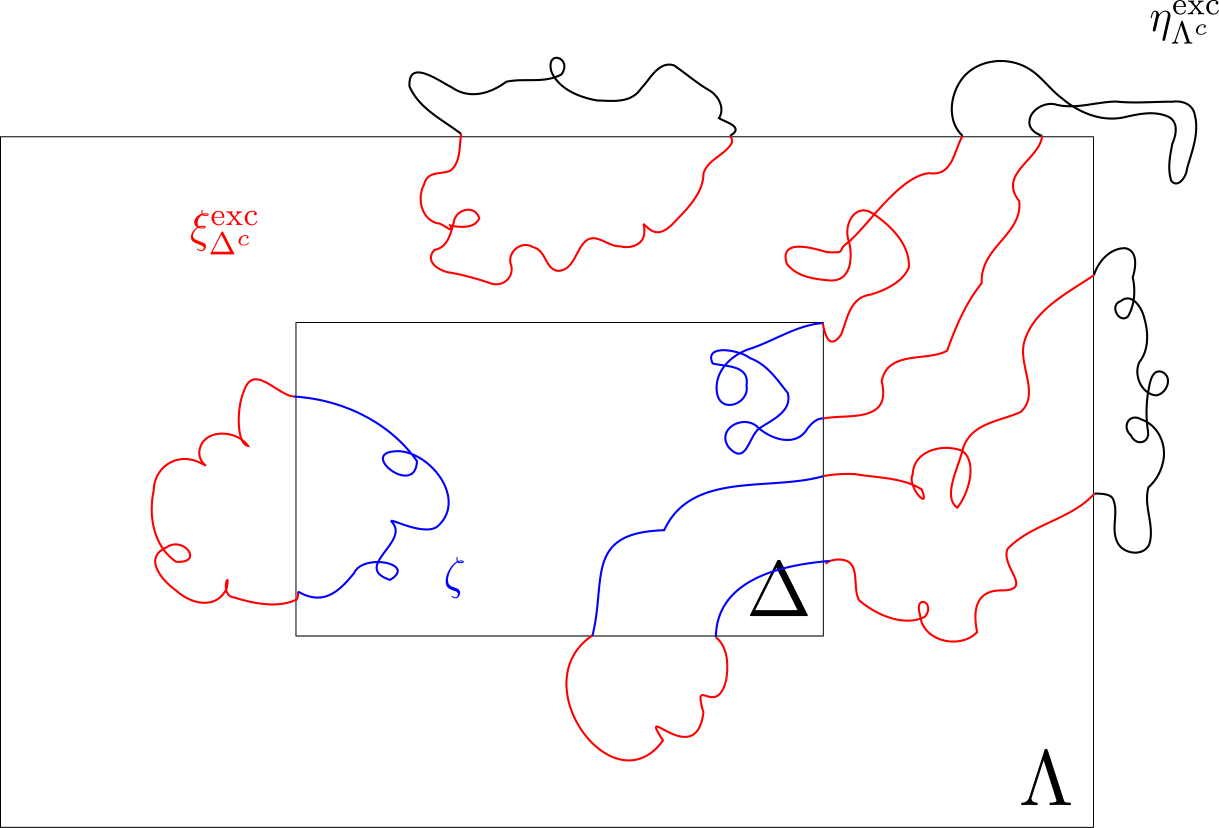}
\caption{An illustration for $\eta^\exc_{\L^c}, \xi^\exc_{\D^c},$ and $\zeta$ in \eqref{eq:Q_consistency}.}
\label{fig:Q_consistency}
\end{figure}
It suffices to show that
\begin{equation}
\int\d\P_\L^\mathrm{exc}(\xi|\eta)f(\xi+\eta^{\mathrm{exc}}_{\L^c})=
\int\d\P_\L^\mathrm{exc}(\xi|\eta)\int\d\P_\D^\mathrm{exc}(\zeta|\xi+\eta^{\mathrm{exc}}_{\L^c})
f(\zeta+\xi^{\mathrm{exc}}_{\D^c}+\eta^{\mathrm{exc}}_{\L^c}),
\end{equation}
which is clear by the Markov property of Brownian excursions, see \cite{sznitman2013scaling} for a similar computation in the context of random Brownian interlacements.
\end{proof}
Moreover, we define the following Gibbs-kernel
\begin{equation}\label{Equation deltaexc}
    \gamma_\L(A|\eta)=\frac{1}{Z^{\mathrm{exc}}_\L\rk{\eta_{\L^c}^{\mathrm{exc}}}}\int\1_A\rk{\eta_{\L^c}^{\mathrm{exc}}+\xi}\ex^{-\beta\H_\L\rk{\eta_{\L^c}^{\mathrm{exc}}+\xi}+\beta\mu \mathrm{N}_\L(\eta_{\L^c}^\exc+\xi)}\d\Q_\L(\xi|\eta)\, ,
\end{equation}
where 
\begin{equation}\label{Equation Zexc}
    Z^{\mathrm{exc}}_\L\rk{\eta_{\L^c}^{\mathrm{exc}}}=\int\ex^{-\beta\H_\L\rk{\eta_{\L^c}^{\mathrm{exc}}+\xi}+\beta\mu \mathrm{N}_\L(\eta_{\L^c}^\exc+\xi)}\d\Q_\L(\xi|\eta)\, .
\end{equation}

\begin{lemma}\label{LemmaConsistentQkernel}
The family $\rk{\gamma_\L}_\L$ is consistent.
\end{lemma}
\begin{proof}
Fix $\D\subset\L$, we want to show that
\begin{equation}
\gamma_\L(\gamma_\D(A|\cdot)|\eta)=\gamma_\L(A|\eta).
\end{equation}
Similar to Lemma \ref{LemmaConsistent},
we assume for simplicity that $\mu=0$. Then
\begin{equation}\label{eq:Htrick_gamma}
\begin{aligned}
\mathrm{Z}_\L^{\mathrm{exc}}&(\eta^{\mathrm{exc}}_{\L^c})\gamma_\L(\gamma_\D(A|\cdot)|\eta)\\
&=\int \d\Q_\L(\xi|\eta)
\gamma(A|\eta^{\mathrm{exc}}_{\L^c}+\xi)\ex^{-\beta\H_\L(\eta_{\L^c}^{\mathrm{exc}}+\xi)}\\
&=\int \d\Q_\L(\xi|\eta)
\int \d\Q_\D(\zeta|\xi+\eta^{\mathrm{exc}}_{\L^c})
\frac{\1_A(\eta^{\mathrm{exc}}_{\L^c}+\xi^{\mathrm{exc}}_{\D^c}+\zeta)}{\mathrm{Z}_\D^{\mathrm{exc}}(\eta^{\mathrm{exc}}_{\L^c}+\xi_{\D^c}^{\mathrm{exc}})}\ex^{-\beta\H_\D(\eta^{\mathrm{exc}}_{\L^c}+\xi^{\mathrm{exc}}_{\D^c}+\zeta)-\beta\H_\L(\eta_{\L^c}^{\mathrm{exc}}+\xi)}, .
\end{aligned}
\end{equation}
By Lemma \ref{lem:boundary_condition} and Lemma \ref{lem:H_trick_exc} we have
\begin{equation}
\H_\D(\eta^{\mathrm{exc}}_{\L^c}+\xi^{\mathrm{exc}}_{\D^c}+\zeta)
+\H_\L(\eta_{\L^c}^{\mathrm{exc}}+\xi)
=\H_\D(\eta^{\mathrm{exc}}_{\L^c}+\xi)
+\H_\L(\eta_{\L^c}^{\mathrm{exc}}+\xi^{\mathrm{exc}}_{\D^c}+\zeta),
\end{equation}
therefore
\eqref{eq:Htrick_gamma} can be further simplified as 
\begin{equation}\label{eq:454}
\begin{aligned}
&\int \d\Q_\L(\xi|\eta)
\int \d\Q_\D(\zeta|\xi+\eta^{\mathrm{exc}}_{\L^c})
\frac{\1_A(\eta^{\mathrm{exc}}_{\L^c}+\xi^{\mathrm{exc}}_{\D^c}+\zeta)}{\mathrm{Z}_\D^{\mathrm{exc}}(\eta^{\mathrm{exc}}_{\L^c}+\xi_{\D^c}^{\mathrm{exc}})}\ex^{-\beta\H_\D(\eta^{\mathrm{exc}}_{\L^c}+\xi)-\beta\H_\L(\eta_{\L^c}^{\mathrm{exc}}+\xi^{\mathrm{exc}}_{\D^c}+\zeta)}.
\end{aligned}
\end{equation}
Apply Lemma \ref{lem:consistency_Q} to  
\begin{equation}
f(\alpha)=
\frac{\ex^{-\beta\H_\D(\alpha)}}{\mathrm{Z}_\D^{\mathrm{exc}}(\alpha^{\mathrm{exc}}_{\D^c})}
\int \d\Q_\D(\beta|\alpha)\1_A(\beta+\alpha^{\mathrm{exc}}_{\D^c})\ex^{-\beta\H_\L(\alpha^{\mathrm{exc}}_{\D^c}+\beta)},
\end{equation}
we have
\begin{equation}
\begin{aligned}
\int &\d\Q_\L(\xi|\eta)\frac{\ex^{-\beta\H_\D(\xi+\eta^{\mathrm{exc}}_{\L^c})}}{\mathrm{Z}_\D^{\mathrm{exc}}(\xi^{\mathrm{exc}}_{\D^c}+\eta^{\mathrm{exc}}_{\L^c})}
\int \d\Q_\D(\beta|\xi+\eta^{\mathrm{exc}}_{\L^c})\1_A(\beta+\xi^{\mathrm{exc}}_{\D^c}+\eta^{\mathrm{exc}}_{\L^c})\ex^{-\beta\H_\L(\xi^{\mathrm{exc}}_{\D^c}+\eta^{\mathrm{exc}}_{\L^c}+\beta)}\\
=&
\int \d\Q_\L(\xi|\eta)\int \d\Q_\D(\zeta|\xi+\eta^{\mathrm{exc}}_{\L^c})
\frac{\ex^{-\beta\H_\D(\zeta+\xi^{\mathrm{exc}}_{\D^c}+\eta^{\mathrm{exc}}_{\L^c})}}{\mathrm{Z}_\D^{\mathrm{exc}}(\zeta^{\mathrm{exc}}_{\D^c}+\xi^{\mathrm{exc}}_{\D^c}+\eta^{\mathrm{exc}}_{\L^c})}\\
&\qquad\qquad\int \d\Q_\D(\beta|\zeta+\xi^{\mathrm{exc}}_{\D^c}+\eta^{\mathrm{exc}}_{\L^c})\1_A(\beta+\zeta^{\mathrm{exc}}_{\D^c}+\xi^{\mathrm{exc}}_{\D^c}+\eta^{\mathrm{exc}}_{\L^c})\ex^{-\beta\H_\L(\zeta^{\mathrm{exc}}_{\D^c}+\xi^{\mathrm{exc}}_{\D^c}+\eta^{\mathrm{exc}}_{\L^c}+\beta)}\\
=&
\int \d\Q_\L(\xi|\eta)\int \d\Q_\D(\zeta|\xi+\eta^{\mathrm{exc}}_{\L^c})
\frac{\ex^{-\beta\H_\D(\zeta+\xi^{\mathrm{exc}}_{\D^c}+\eta^{\mathrm{exc}}_{\L^c})}}{\mathrm{Z}_\D^{\mathrm{exc}}(\xi^{\mathrm{exc}}_{\D^c}+\eta^{\mathrm{exc}}_{\L^c})}\\
&\qquad\qquad\int \d\Q_\D(\beta|\xi+\eta^{\mathrm{exc}}_{\L^c})\1_A(\beta+\xi^{\mathrm{exc}}_{\D^c}+\eta^{\mathrm{exc}}_{\L^c})\ex^{-\beta\H_\L(\xi^{\mathrm{exc}}_{\D^c}+\eta^{\mathrm{exc}}_{\L^c}+\beta)}\\
=&\int \d\Q_\L(\xi|\eta)
\int \d\Q_\D(\beta|\xi+\eta^{\mathrm{exc}}_{\L^c})\1_A(\beta+\xi^{\mathrm{exc}}_{\D^c}+\eta^{\mathrm{exc}}_{\L^c})\ex^{-\beta\H_\L(\xi^{\mathrm{exc}}_{\D^c}+\eta^{\mathrm{exc}}_{\L^c}+\beta)}.
\end{aligned}
\end{equation}
Therefore, \eqref{eq:454} is equal to
\begin{equation}\label{eq:457}
\int \d\Q_\L(\xi|\eta)
\int \d\Q_\D(\beta|\xi+\eta^{\mathrm{exc}}_{\L^c})\1_A(\beta+\xi^{\mathrm{exc}}_{\D^c}+\eta^{\mathrm{exc}}_{\L^c})\ex^{-\beta\H_\L(\xi^{\mathrm{exc}}_{\D^c}+\eta^{\mathrm{exc}}_{\L^c}+\beta)}.
\end{equation}
Use Lemma \ref{lem:consistency_Q} again for
\begin{equation}
f(\alpha)=\1_A(\alpha)\ex^{-\beta\H_\L(\alpha)},
\end{equation}
we further reduce \eqref{eq:457} to
\begin{equation}
\int \d\Q_\L(\xi|\eta)
\1_A(\xi+\eta^{\mathrm{exc}}_{\L^c})\ex^{-\beta\H_\L(\xi+\eta^{\mathrm{exc}}_{\L^c})}
=\mathrm{Z}_\L^{\mathrm{exc}}(\eta^{\mathrm{exc}}_{\L^c})\gamma_\L(A|\eta).
\end{equation}
In conclusion, we have
\begin{equation}
\mathrm{Z}_\L^{\mathrm{exc}}(\eta^{\mathrm{exc}}_{\L^c})\gamma_\L(\gamma_\D(A|\cdot)|\eta)=\mathrm{Z}_\L^{\mathrm{exc}}(\eta^{\mathrm{exc}}_{\L^c})\gamma_\L(A|\eta).
\end{equation}
\end{proof}
Next, we show that the excursion partition function is almost surely finite with respect to ${\gfrak}_n$
\begin{lemma}\label{Lem:ExcPartFiniteFinite}
For all $n\ge 1$ and all $\Delta\subset\L_n$, we have that
\begin{equation}
    \gfrak_n\rk{\mathrm{Z}_\D^{\mathrm{exc}}=0}=0\, .
\end{equation}
\end{lemma}
\begin{proof}
     We follow \cite{dereudre2012existence}. Note that by construction $\tilde{\gfrak}_n\rk{\H_{\L_n}=\infty}=0$. In Equation \eqref{eq:Htrick_gamma} we plug in $\L=\L_n$, $\Delta=\Delta$ and $\eta=\emptyset$ to obtain
    \begin{equation}
        \ex^{-\H_{\L_n}(\xi)}\ex^{-\H_\D(\xi_{\D^c}^{\mathrm{exc}}+\zeta)}=\ex^{-\H_{\L_n}(\xi_{\D^c}^{\mathrm{exc}}+\zeta)}\ex^{-\H_\D(\xi)}\, .
    \end{equation}
    By integrating over $\zeta\subset\D$, we obtain
    \begin{equation}
        \ex^{-\H_{\L_n}(\xi)}\mathrm{Z}_\D^\mathrm{exc}\rk{\xi}=\mathrm{Z}_{\L_n,\D}^\mathrm{exc}(\xi)\ex^{-\H_\D(\xi)}\, ,
    \end{equation}
    where $\mathrm{Z}_{\L_n,\D}(\xi)=\int \ex^{-\H_{\L_n}\rk{\xi_{\D}^{\dir,c}+\zeta}}\d \Q_\D\rk{\zeta|\xi_{\D^c}^{\mathrm{exc}}}$. Note that $\d\P_{\L_n}^\dir(\xi)=\d\Q_{\L_n}(\xi|\emptyset)$ and that hence for any test function $F$
    \begin{equation}
        \int\d\P_{\L_n}^\dir(\xi)F(\xi)= \int\d\P_{\L_n}^\dir(\xi)\int \d\Q_{\L_n}\rk{\zeta|\xi_{\D^c}^\mathrm{exc}}F\rk{\xi_{\D^c}^{\mathrm{exc}}+\zeta}\, .
    \end{equation}
    Hence
    \begin{equation}
    \begin{split}
        \tilde{\gfrak}_n\rk{\mathrm{Z}_\D^{\mathrm{exc}}=0}&= \tilde{\gfrak}_n\rk{\mathrm{Z}_{\L_n,\D}^{\mathrm{exc}}(\xi)=0}=\frac{1}{Z_{\L_n}(\emptyset)}\int\ex^{-\H_{\L_n}\rk{\xi}}\1\rk{\mathrm{Z}_{\L_n,\D}(\xi)=0}\d \P_{\L_n}^\dir(\xi)\\
        &=\frac{1}{\mathrm{Z}_{\L_n}(\emptyset)}\int \mathrm{Z}_{\L_n,\D}^{\mathrm{exc}}(\xi_{\D^c}^{\mathrm{exc}})\1\rk{\mathrm{Z}_{\L_n,\D}^{\mathrm{exc}}(\xi_{\D^c}^{\mathrm{exc}})=0}\d\P_{\L_n}^\dir(\xi)=0\, .
    \end{split}
    \end{equation}
    Since $\tilde{\gfrak}_n\rk{\mathrm{Z}_\D^{\mathrm{exc}}=0}=0$, we immediately get that ${\gfrak}_n\rk{\mathrm{Z}_\D^{\mathrm{exc}}=0}=0$. This concludes the proof.
\end{proof}

\subsection{Entropy}\label{sec:4.4}
In this section we introduce the specific entropy (also called mean entropy per site) and prove an upper bound for the entropy of $\gfrak_n$. By the compactness of level sets, this gives us an accumulation point of $\rk{\gfrak_n}_n$.
\begin{definition}
For two probability measures, $P\ll Q$, we define the relative entropy as
\begin{align}
h(P|Q)=\int \d P\log\frac{\d P}{\d Q}.
\end{align}
\end{definition}
Moreover, for a probability measure $P$ on $\Omega$, write $P_\L$ for the law induced by the map $\eta\mapsto\eta_\L$. We define the average relative entropy on $\L$ as
\begin{equation}
    I_\L(P)=\frac{h\rk{P_\L|\P_\L}}{\abs{\L}}=\frac{1}{|\L|}\int \d P_\L \log\frac{\d P_\L}{\d\P_\L} 
    \, .
\end{equation}
We set $I$ the \textit{specific entropy}
\begin{equation}\label{Equation I}
    I(P)=\lim_{n\to\infty}I_{\L_n}(P)\quad\text{if the limit exists,}
\end{equation}
where $\L_n$ is the centred cube of side-length $n$ in $\mathbb R^d$.
\begin{lemma}\label{lementrpyExists}
For $P$ translation invariant with respect to translations in $\Z^d$, $I(P)$ is well-defined and equals
\begin{equation}
    I(P)=\sup_{n\in\N}I_{\L_n}(P)\, .
\end{equation}
\end{lemma}
This follows from \cite[Chapter 15]{georgii1988gibbs}.

\begin{lemma}\label{LemmaEntropyEstimates}
For every $n>0$, $I(\gfrak_n)$ is well-defined. Furthermore if Assumption \ref{refassump} is satisfied, then there exists a $C_1>0$ such that
\begin{equation}
    I(\gfrak_n)\le C_1\quad\textnormal{for all}\quad n\ge 0.
\end{equation}
\end{lemma}
\begin{proof}
{It suffices to prove for $I(\hat\gfrak_n)$, since by \cite[Chapter 15]{georgii1988gibbs}
\begin{equation}\label{eq:temp484}
    I(\gfrak_n)=\frac{1}{n^d}\sum_{x\in\L_n\cap\Z^d}I(\hat\gfrak_n\circ\tau_x)=I(\hat\gfrak_n)\, .
\end{equation}}

By definition, for disjoint $\L$ and $\L'$,
\begin{equation}
\frac{\d\delta^\dir_{\L}(\cdot|\emptyset)\otimes\delta^\dir_{\L'}(\cdot|\emptyset)}{\d\P^\dir_{\L\cup\L'}}
=\frac{\d\delta^\dir_{\L}(\cdot|\emptyset)}{\d\P^\dir_{\L}}\cdot\frac{\d\delta^\dir_{\L'}(\cdot|\emptyset)}{\d\P^\dir_{\L'}},
\text{ and }
\frac{\d\P^\dir_{\L\cup\L'}}{\d\P_{\L\cup\L'}}=\frac{\d\P^\dir_{\L}}{\d\P_{\L}}\cdot\frac{\d\P^\dir_{\L'}}{\d\P_{\L'}}.
\end{equation}
So if we take any $n,k\in\mathbb N$, and let $m=(2k+1)n$, then
{\small\begin{equation}
\begin{aligned}
&I_{\L_m}(\hat\gfrak_n)\\
=&\frac{1}{m^d}h(\hat\gfrak_n|\P_{\L_m})\\
=&\frac{1}{m^d}\int\log
\frac{\d\bigotimes_{-k\le x\le k}\delta^\dir_{\L_n+xn}(\cdot|\emptyset)}
{\d\P_{\L_m}}
\d\bigotimes_{-k\le x\le k}\delta^\dir_{\L_n+xn}(\cdot|\emptyset)\\
=&\frac{1}{m^d}\int\log
\frac{\d\bigotimes_{-k\le x\le k}\delta^\dir_{\L_n+xn}(\cdot|\emptyset)}
{\d\P_{\L_m}^\dir}
\d\bigotimes_{-k\le x\le k}\delta^\dir_{\L_n+xn}(\cdot|\emptyset)
+\frac{1}{m^d}\int\log
\frac{\d\P_{\L_m}^\dir}
{\d\P_{\L_m}}
d\bigotimes_{-k\le x\le k}\delta^\dir_{\L_n+xn}(\cdot|\emptyset)\\
=&\frac{1}{n^d}\int\log
\frac{\d\delta_{\L_n}(\cdot|\emptyset)}
{\d\P_{\L_m}^\dir}
\d\delta^\dir_{\L_n}(\cdot|\emptyset)
+\frac{1}{n^d}\int\log
\frac{\d\P_{\L_n}^\dir}
{\d\P_{\L_n}}
\d\delta^\dir_{\L_n}(\cdot|\emptyset)\\
=&I_{\L_n}(\tilde\gfrak_n).
\end{aligned}
\end{equation}}
Take $m\rightarrow\infty$, we then have
\begin{equation}
I(\hat\gfrak_n)=I_{\L_n}(\tilde\gfrak_n).
\end{equation}

Furthermore, notice that
\begin{equation}
    \frac{\d\P_\L^\dir}{\d\P_\L}(\eta)=\ex^{\mathrm{M}_\L[\omega\cap \L^c\neq \emptyset]}\1\{\eta=\eta_\L^\dir\}\, ,
\end{equation}
and we deduce that
\begin{equation}\label{eq:three_terms}
\begin{aligned}
    {n^d}I(\hat\gfrak_n)=n^d&{I_{\L_n}(\tilde\gfrak_n)}=h(\tilde\gfrak_n|\P_{\L_n})=h(\tilde\gfrak_n|\P_{\L_n}^\dir)+\tilde\gfrak_n\ek{\log\frac{\d\P_{\L_n}^\dir}{\d\P_{\L_n}}}\\
    &=-\log \mathrm{Z}_{\L_n}^\dir(\emptyset)-\int \beta\rk{\H_{\L_n}(\xi)-\mu \mathrm{N}_{\L_n}(\xi)}\delta^\dir_{\L_n}(\d \xi|\emptyset)+\mathrm{M}^\dir_{\L_n}\ek{\omega\cap\L_n^c\neq\emptyset}\, .
\end{aligned}
\end{equation}
For the first term, note that there exists $c>0$ such that
\begin{equation}\label{eq:empty_bound}
\begin{aligned}
    \mathrm{Z}_{\L_n}^\dir(\emptyset)=\int\ex^{-\beta\H_{\L_n}(\xi)+\beta\mu \mathrm{N}_{\L_n}(\xi)}\d \P^\dir_{\L_n}(\xi)
    &\ge \P^\dir_{\L_n}(\xi=\emptyset)\\
    &=\ex^{-\mathrm{M}_{\L_n}^\dir[\1]}\\
    &= \ex^{-\int_{\L_n}\d x\sum_{j}\frac{1}{j}p_{\beta j}\hk{\L_n}(x,x)}
    \ge \ex^{-c\abs{\L_n}}.
\end{aligned}
\end{equation}
For the second term, by Lemma \ref{lem:FKG},
\begin{equation}
\begin{split}
  -\int  (\beta\H_{\L_n}(\xi)-\mu \mathrm{N}_{\L_n}(\xi))\delta^\dir_{\L_n}(\d \xi|\emptyset)    \le \int \beta\abs{\mu} \mathrm{N}_{\L_n}(\xi)\delta^\dir_{\L_n}(\d \xi|\emptyset)\le \beta|\mu|\E^\H_{\L_n}[\mathrm{N}_{\L_n}(\xi)]\,,
\end{split}
\end{equation}
and  by Campbell's formula
\begin{equation}
\begin{split}
  \E^\H_{\L_n}[\mathrm{N}_{\L_n}(\xi)]=M^\H_{\L_n}\ek{\ell(\omega)}\le \int_{\L_n}\d x\sum_{j\ge 1}{\ex^{\beta \mu j}}\ex^{-c_\Phi j}p_{\beta j}(0)=\Ocal\rk{n^d}\, .
\end{split}
\end{equation}
To estimate the third term $\mathrm{M}^\dir_{\L_n}\ek{\omega\cap\L_n^c\neq\emptyset}$, note that
\begin{equation}
\P_{x,x}^t(\omega\cap \L^c\ne\emptyset)
\le d\cdot \P_{x,x}^t\left(\sup_{0\le s\le t}|\omega_1(s)-x|\ge \frac{\dist(x,\L^c)}{\sqrt d}\right),
\end{equation}
where $\omega_1$ is the first coordinate of $\omega$,
so it is clear that when $R\rightarrow \infty$, 
uniformly for all $x$ such that $\dist(x,\L^c)\ge R$,
\begin{equation}
\sum_{j\ge 1}\frac{1}{j}\mathbb P_{x,x}^{\beta j}(\omega\cap\L^c\ne\emptyset)\rightarrow 0\, .
\end{equation}
So for any $\epsilon>0$, we can find large enough $R>0$, such that
\begin{equation}
\begin{aligned}
    \frac{1}{|\L|}&\mathrm{M}^\dir_\L[\omega\cap\L^c\neq \emptyset]\\
&=   \frac{1}{|\L|}\int_\L \d x\sum_{j\ge 1}\frac{1}{j}\mathbb P_{x,x}^{\beta j}(\omega\cap\L^c\ne\emptyset)\\
&\le
\epsilon+\frac{\int_\L \1(\dist(x,\L^c)<R)\d x}{|\L|}\cdot\sum_{j\ge 1}\frac{1}{j}\mathbb E_{x,x}^{\beta j}[\1].
\end{aligned}
\end{equation}
Since
\begin{equation}
\lim_{\L\uparrow\mathbb R^d}
\frac{\int_\L \1(\dist(x,\L^c)<R)\d x}{|\L|}=0,
\end{equation}
we conclude that
\begin{equation}
\lim_{\L\uparrow\mathbb R^d}\frac{1}{|\L|}\mathrm{M}^\dir_\L[\omega\cap\L^c\neq \emptyset]=0.
\end{equation}

Combining the estimates above for the three terms in \eqref{eq:three_terms}, we have
\begin{equation}
    n^dI_{\L_n}(\tilde\gfrak_n)=\Ocal(n^d)\, .
\end{equation}
This concludes the proof.
\end{proof}
The following corollary is standard, see \cite[Proposition 2.6]{georgii1993large}.
\begin{cor}\label{CorollaryCompact}
Under the assumption of Lemma \ref{LemmaEntropyEstimates}, there exists a probability measure $\gfrak$ on $\Omega$ such that for every bounded function $F$, depending only on $\eta_\L$ for some compact $\L$,
\begin{equation}
    \lim_{n\to\infty}\gfrak_n\ek{F}=\gfrak\ek{F}\, .
\end{equation}
\end{cor}
\begin{proof}
By the above mentioned reference, the level sets $\gk{P\colon I(P)\le c}$ (for $c<\infty$) are compact in the topology of continuous and bounded cylinder functions.
\end{proof}
\begin{remark}\label{remarksuperstable}
If $\Phi$ is no longer non-negative but simply \textit{superstable} (see \cite{ruelle1999statistical} for a definition), the above result remains valid. As we lack a proof that $\gfrak$ is Gibbs for all superstable potentials, we chose to defer this case for future investigations.
\end{remark}
Next, we seek to strengthen the topology in which $\rk{\gfrak_n}_n$ converges. This step will be crucial in establishing the Gibbs property.
\subsection{Topology}\label{sec:4.5}
For a translation invariant function $\psi\colon\Gamma\to\R$, we abbreviate
\begin{equation}
   \mathrm{N}_\L^\psi(\eta)=\sum_{\omega\in\eta_\L}\psi(\omega)\, .
\end{equation}
We say that a function $F\colon\Omega\to\R$ is \textit{local}, if there exists $\D$ bounded, such that $F$ depends only on $\eta_\L$.
\begin{definition}
We say that a local function $F$ (with support $\L$) is \textnormal{tame} (or $\psi$-tame), if there exists a $C=C_F>0$ such that
\begin{equation}
    \abs{F(\eta)}\le C\rk{1+\mathrm{N}_\L^\psi(\eta)}\, .
\end{equation}
A sequence of $(Q_n)_n$ probability measures converges to $Q$ \textnormal{in the topology of local convergence} if for every local and tame function $F$,
\begin{equation}
    \lim_n Q_n\ek{F}=Q[F]\, .
\end{equation}
We denote this topology by $\LC=\LC_\psi$.
\end{definition}
For example, $\psi\equiv 0$ gives the topology of local and bounded functions. $\psi(\omega)=\ell(\omega)$ gives the topology of local functions which grow at most proportionally with the particle number $\mathrm{N}$.

In order to examine in what topology $\rk{\gfrak_n}_n$ may converge, we need to bound the expectation of certain $\psi$ functions.

\begin{lemma}\label{LemmmaLoopLengths}
Suppose $f\colon\N\to (0,\infty)$ such that
\begin{equation}
    \sum_{j\ge 1}\frac{\ex^{j(-c_\Phi+\beta\mu)}}{j}f(j)<\infty\, .
\end{equation}
Then
\begin{equation}
    \gfrak_n\ek{\sum_{\omega\in \eta}f(\ell(\omega))}=\Ocal\rk{n^d}\, .
\end{equation}
\end{lemma}
\begin{proof}
Since $\mathbb P^\dir_{\L_n}$ is a Poisson point process, for any function $G$,
\begin{equation}\label{eq:60}
\mathbb E_{\L_n}^\dir[G]=\ex^{-\mathrm{M}_{\L_n}^\dir[\Gamma]}
\sum_{k\ge 0}\frac{1}{k!}\rk{\mathrm{M}_{\L_n}^\dir}^{\otimes k}[G]
\, ,
\end{equation}
see \cite[Chapter 3]{last2017lectures}. Here, $G\big((\omega_1,\ldots,\omega_k)\big)$ is understood as $G\rk{\delta_{\omega_1}+\ldots+\delta_{\omega_k}}$.
In particular,
\begin{equation}
\mathrm{Z}_{\L_n}^\dir(\emptyset)={\mathbb E_{\L_n}^\dir[\ex^{-\beta \H_{\L_n}+\beta\mu \mathrm{N}_{\L_n}}]}=\ex^{-\mathrm{M}_{\L_n}^\dir[\Gamma]}\sum_{k\ge 0}\frac{1}{k!}\rk{\mathrm{M}_{\L_n}^\dir}^{\otimes k}[\ex^{-\beta\H_{\L_n}+\beta\mu \mathrm{N}_{\L_n}}].
\end{equation}
For the rest of this proof, abbreviate $\mathrm{M}_{\L_n}^\dir$ by $\mathrm{M}$.

By definition of $\gfrak_n$, for any function $F$
\begin{equation}\label{eq:59}
\gfrak_n\left[F\right]
=
\frac{1}{\mathrm{Z}^\dir_{\L_n}(\emptyset)}\mathbb E_{\L_n}^\dir\left[\ex^{-\beta \H_{\L_n}+\beta\mu \mathrm{N}_{\L_n}}F\right].
\end{equation}
Then,
{\small\begin{equation}
    \gfrak_n\left[\sum_{j\ge 1}f(j)\#\gk{\omega\colon \ell(\omega)=j}\right]=\frac{1}{\mathrm{Z}^\dir_{\L_n}(\emptyset)}\ex^{-\mathrm{M}[\Gamma]}\sum_{j\ge 1}\sum_{k\ge 1}\frac{f(j)}{k!}\mathrm{M}^{\otimes k}\ek{\#\gk{\omega\colon \ell(\omega)=j}\ex^{-\beta\H_{\L_n}+\beta\mu \mathrm{N}_{\L_n}}}.
\end{equation}}
We also have that
\begin{equation}
\begin{aligned}
    \mathrm{M}^{\otimes k}\ek{\#\gk{\omega\colon \ell(\omega)=j}\ex^{- \beta\H_{\L_n}+\beta\mu \mathrm{N}_{\L_n}}}
    &=\mathrm{M}^{\otimes k}\ek{\ex^{-\beta\H_{\L_n}+\beta\mu \mathrm{N}_{\L_n}}\sum_{i=1}^k \1\{\ell(\omega_i)=j\}}
    \\
    &=kM^{\otimes k}\ek{\1\{\ell(\omega_1)=j\}\ex^{-\beta\H_{\L_n}+\beta\mu \mathrm{N}_{\L_n}}}\, .
\end{aligned}
\end{equation}
Since we have non-negative interaction, we get that
\begin{equation}
\begin{aligned}
    \mathrm{M}^{\otimes k}&\ek{\1\{\ell(\omega_1)=j\}\ex^{-\beta\H_{\L_n}+\beta\mu \mathrm{N}_{\L_n}}}\\
    \le &\mathrm{M}^{\otimes (k-1)}\ek{\ex^{(-\beta\H_{\L_n}+\beta\mu \mathrm{N}_{\L_n})(\omega_1,\ldots,\omega_{k-1})}}\mathrm{M}\ek{\1\{\ell(\omega_1)=j\}\ex^{(-\beta\H_{\L_n}+\beta\mu \mathrm{N}_{\L_n})(\omega_1)}}\, .
\end{aligned}
\end{equation}
Note that we have by using Lemma \ref{LemmaSingleLoopDecay} in the last equality
\begin{equation}
    \mathrm{M}\ek{\1\{\ell(\omega_1)=j\}\ex^{(-\beta\H_{\L_n}+\beta\mu \mathrm{N}_{\L_n})(\omega_1)}}\le \frac{\abs{\L_n}\ex^{\beta \mu j}}{j}\E_{0,0}\hc{\L}{\beta j}\left[\ex^{-\beta \H_{\L_n}}\right]=n^dj^{-1} \Ocal\left(\ex^{j(-c_\Phi +\beta\mu)}\right)\, .
\end{equation}
Summarising the above, we get that 
\begin{equation}
\begin{aligned}
\gfrak_n&\left[\sum_{j\ge 1}f(j)\#\gk{\omega\colon \ell(\omega)=j}\right]\\
&\le\frac{C}{\mathrm{Z}_{\L_n}^\dir(\emptyset)}\ex^{-\mathrm{M}[\Gamma]}\sum_{j\ge 1}\sum_{k\ge 1}\frac{f(j)}{k!}\frac{kn^d\ex^{j(-c_\Phi +\beta\mu)}}{j}\mathrm{M}^{\otimes (k-1)}[\ex^{-\beta\H_{\L_n}+\beta\mu \mathrm{N}_{\L_n}}]\\
&=
Cn^d \sum_{j\ge 1}\frac{f(j)\ex^{j(-c_\Phi +\beta\mu)}}{j}\cdot
\frac{\ex^{-\mathrm{M}[\Gamma]}}{\mathrm{Z}^\dir_{\L_n}(\emptyset)}\sum_{k\ge 1}\frac{\mathrm{M}^{\otimes (k-1)}[\ex^{-\beta\H+\beta\mu N}]}{(k-1)!}\\
&=Cn^d\sum_{j\ge 1} f(j)\frac{\ex^{j(-c_\Phi +\beta\mu)}}{j}.
\end{aligned}
\end{equation}
\end{proof}
Recall that a function $F$ is $\psi$-tame if there exist $C>0$ such that $\abs{F}\le C\rk{1+N_\L^\psi}$. By Corollary \ref{CorollaryCompact}, we know that $\gfrak_n$  always converges to $\gfrak$ in the $\LC_\psi$ topology for $\psi(\omega)=0$. We now try to lift it to larger $\psi$.
\begin{proposition}\label{PropTopology}
Let $U$ be the centred unit cube. For every function $\psi\colon\Gamma\to(0,\infty)$ so that
\begin{equation}
    \gfrak\ek{\mathrm{N}_U^\psi}<\infty\, ,
\end{equation}
if there exists a function $f\colon\N\to(0,\infty)$ with the following conditions:
\begin{itemize}
\item The function $f$ satisfies
\begin{equation}
    \sum_{j\ge 1}\frac{\ex^{j(-c_\Phi+\beta\mu)}}{j}f(j)<\infty\, .
\end{equation}
\item For every $a>0$, there exists $C_a>0$ so that for every $\D\subseteq\mathbb R^d$ compact
\begin{equation}\label{Equationcodnitiontop}
    \int_{\D}\sum_{j\ge 1}\frac{\ex^{- f(j)}}{j}\E_{x,x}^{\beta j}\ek{\ex^{a\psi}}\d x\le C_a\abs{\D}\, .
\end{equation}
\end{itemize}
Then
\begin{equation}
    \gfrak_n\to\gfrak\quad\text{in}\quad \LC_\psi\, .
\end{equation}
\end{proposition}
\begin{proof}
In fact, it suffices to show for any compact set $\Delta$, that
\begin{equation}\label{Equationbtoinfty}
    \lim_{b\to\infty}\gfrak_n\ek{\mathrm{N}_\D^\psi\1\gk{\mathrm{N}_\D^\psi>b}}=0\quad\emph{uniformly}\text{ in }n\in\N\, .
\end{equation}
Indeed, we can estimate
{\small\begin{equation}\label{eq:g_n_to_g}
\begin{aligned}
    |&\gfrak_{n}[F]-\gfrak[F]|\\
&\le |\gfrak_{n}[F]-\gfrak_{n}[F\1\{\mathrm{N}_\D^\psi\le b\}]|+|\gfrak_{n}[F\1\{\mathrm{N}_\D^\psi\le b\}]-\gfrak[F\1\{\mathrm{N}_\D^\psi\le b\}]|+|\gfrak[F\1\{\mathrm{N}_\D^\psi\le b\}]-\gfrak[F]|\, .
\end{aligned}
\end{equation}}
For $\psi$-tame $F$, this shows the equivalence of Proposition \ref{PropTopology} and Equation \eqref{Equationbtoinfty}.

A vital tool for us will be the entropy inequality: for any two probability measures $\mu$ and $\nu$ on some Polish space $\Sigma$
\begin{equation}
    h(\mu|\nu)\ge \mu[F]-\log\nu\ek{\ex^{F}}\, ,
\end{equation}
for any $F$ measurable and bounded, see \cite[Lemma 6.2.13]{dembo2009large}. 
We will use this formula with $\mu=\gfrak_n$. However, we cannot choose $\nu=\P_{\L_n}$. This is because $\mathrm{N}^\psi_\L$ has worse integrability properties under $\P_{\L_n}$ than under $\gfrak_n$. For this, we define
\begin{equation}\label{EquationIntermediateMeause}
    \tilde{\mathrm{M}}_{\L_n}=\int_{\L_n}\sum_{j\ge 1}\frac{\ex^{- f(j)}}{j}\P_{x,x}^{\beta j}\,\,\d x\, .
\end{equation}
Let $\tilde\P_{\L_n}$ be the Poisson point process with intensity measure $ \tilde{\mathrm{M}}_{\L_n}$.

Note that
\begin{equation}
    \frac{\d \P_{\L_n}}{\d \tilde\P_{\L_n}}(\eta)=\ex^{ \tilde{\mathrm{M}}_{\L_n}\ek{\Gamma}-\mathrm{M}_{\L_n}\ek{\Gamma}}\exp\rk{\sum_{j\ge 1}f(j)n_j(\eta)}\, ,
\end{equation}
where $n_j(\eta)$ is short for $\#\{\omega\in\eta:\ell(\omega)=j\}$,
and thus
\begin{equation}
\begin{aligned}
    h&\rk{\gfrak_n|\tilde\P_{\L_n}}\\
    &=h\rk{\gfrak_n|\P_{\L_n}}+\gfrak_n\ek{\log  \frac{\d \P_{\L_n}}{\d \tilde\P_{\L_n}}}=h\rk{\gfrak_n|\P_{\L_n}}+ \tilde{\mathrm{M}}_{\L_n}\ek{\Gamma}-\mathrm{M}_{\L_n}\ek{\Gamma}+\gfrak_n\ek{\sum_{j\ge 1}f(j)n_j}\, .
\end{aligned}
\end{equation}
As $\tilde\P_{\L_n}$ is a Poisson point process, for every $m\in\N$ the limit
\begin{equation}
    \lim_{n\to\infty}\frac{h\rk{\gfrak_m|\tilde\P_{\L_n}}}{n^d}\, ,
\end{equation}
exists and is equal to $\sup_{n\to\infty}\frac{h\rk{\gfrak_m|\tilde\P_{\L_n}}}{n^d}$, see \cite[Chapter 15]{georgii1988gibbs}. We denote this value by $\tilde I(\gfrak_m)$. By Lemma \ref{LemmmaLoopLengths} and Lemma \ref{LemmaEntropyEstimates}, and the fact that $ \tilde{\mathrm{M}}_{\L_n}\ek{\Gamma}$ and $\mathrm{M}_{\L_n}\ek{\Gamma}$ are both of order $n^d$, we have
\begin{equation}
    \sup_{m\in\N} \tilde I(\gfrak_m)\le C\, ,
\end{equation}
for some $C>0$.
For any function $F$ which is $\Fcal_{\L_m}$-measurable, we have by the previous equality, the following inequality holds \textit{for all} $n\in\N$
\begin{equation}
    \gfrak_n[F]-\log\tilde\P_{\L_n}\ek{\ex^F}\le h\rk{\rk{\gfrak_n}_{\L_m}|\rk{\tilde\P_{\L_n}}_{\L_m}}\le  C m^d\, .
\end{equation}
Fix $\e>0$ and $\D\subset\L_m$ for some $m>0$ (also fixed). Choosing $F=a \mathrm{N}_\D^\psi\1\gk{\mathrm{N}_\D^\psi>b}$, with $a=2Cm^d/\e$, then for every $n\in\N$,
\begin{equation}
    \gfrak_n\left[\mathrm{N}_\D^\psi\1\gk{\mathrm{N}_\D^\psi>b}\right]\le \frac{C m^d}{a}+\frac{\log\tilde\P_{\L_n}\ek{\ex^{a\mathrm{N}_\D^\psi\1\gk{\mathrm{N}_\D^\psi>b}}}}{a}\, .
\end{equation}
By Campbell's formula,
\begin{equation}
    \log\tilde\P_{\L_n}\ek{\ex^{a\mathrm{N}_\D^\psi}}= \int_{\L_m}\sum_{j\ge 1}\frac{\ex^{- f(j)}}{j}\P_{x,x}^{\beta j}\ek{\ex^{a\psi}-1}\d x\, ,
\end{equation}
since $\mathrm{N}_\D^\psi$ only depends on values inside $\L_m$. Thus by \eqref{Equationcodnitiontop}, we can choose $b>0$ such that
\begin{equation}
    \log\tilde\P_{\L_n}\ek{\ex^{a\mathrm{N}_\D^\psi\1\gk{\mathrm{N}_\D^\psi>b}}}\le C m^d\, ,
\end{equation}
and therefore
\begin{equation}
\gfrak_n\left[\mathrm{N}_\D^\psi\1\gk{\mathrm{N}_\D^\psi>b}\right]\le\epsilon.
\end{equation}
This concludes the proof of \eqref{Equationbtoinfty}.
\end{proof}
\begin{remark}\label{Rmk:Equationbtoinfty_hat}
Note that $\sup_{n\in\mathbb N}I(\hat \gfrak_n\circ\tau_x)<\infty$, we may replace $\gfrak_n$ by $\hat \gfrak_n\circ\tau_x$ in the arguments above and deduce that  
\begin{equation}
    \hat\gfrak_n\circ\tau_x\ek{\mathrm{N}_\D^\psi}<\infty
\end{equation}
uniformly in all $n\in\mathbb N$ and $x\in\mathbb Z^d$.
\end{remark}
We now use the previous proposition to prove integrability properties.
\begin{lemma}\label{LemmaConvergenceParticlePower}
For $c<c_\Phi$ and any $\D\subset\R^d$ bounded
\begin{equation}
   \gfrak\ek{\ex^{c \mathrm{N}_\D}}<\infty\, .
\end{equation}
Hence, we have $\gfrak_n\to\gfrak$ in $\LC_\psi$ with $\psi(\omega)=\ell(\omega)^p$ for every $p>0$.
\end{lemma}
\begin{proof}
Suppose without loss of generality that $\D=U=[0,1]^d$. Note that by Corollary \ref{CorollaryCompact},

\begin{equation}
    \gfrak\rk{\mathrm{N}_U>k}=\lim_{n\to\infty}\gfrak_n\rk{\mathrm{N}_U>k}\, .
\end{equation}
Note that
\begin{equation}
    \gfrak_n\rk{\mathrm{N}_U>k}\le \ex^{-kc_\Phi}\gfrak_n\ek{\ex^{c_\Phi \mathrm{N}_U }}\, .
\end{equation}
By Lemma \ref{lem:FKG}, for every $n\in\mathbb N$,
\begin{equation}\label{eq:finite_exponential_N}
    \gfrak_n\ek{\ex^{c_\Phi\mathrm{N}_U}}\le \E_{U}^\H\ek{\ex^{c_\Phi\mathrm{N}_U}}=\exp\rk{M^\H_U\ek{\ex^{c_\Phi\ell(\omega)}}}<\infty.
\end{equation}
Hence
\begin{equation}
    \gfrak\rk{\mathrm{N}_U>k}=\Ocal\rk{\ex^{-kc_\Phi}}\, ,
\end{equation}
and the first claim follows.

The second claim follows from the fact that $\gfrak[\mathrm{N}_U^p]$ is finite for any $p>0$ and Proposition \ref{PropTopology} (use for instance $f(j)=j^{p+1}$).
\end{proof}
Recall the definition of $S_\D$ from Equation \eqref{definitinvol}.
\begin{lemma}\label{lemmacontroldiam}
For all $\alpha>0$ and for $\D$ compact
\begin{equation}
    \lim_{n\to\infty}\gfrak_n\ek{\mathrm{S}_\D^\alpha}=\gfrak\ek{\mathrm{S}_\D^\alpha}<\infty\, .
\end{equation}
\end{lemma}
\begin{proof}

Note that for $t>0$, we can find $C>0$ such that
{\small\begin{equation}\label{eq:diameter_use_bridge_to_walk}
\begin{aligned}
        \P_{0,0}^{\beta j}\rk{\diam(\omega)>t}
        &=\P_{0,0}^{1}\rk{\diam(\omega)>t(\beta j)^{-1/2}}\\
        &\le\P_{0,0}^{1}\rk{\diam(\omega[0,1/2])>t(\beta j)^{-1/2}/2}+
        \P_{0,0}^{1}\rk{\diam(\omega[1/2,1])>t(\beta j)^{-1/2}/2}\\
        &\le 2\P_0\rk{\diam(\omega{[0,1/2]})>t(\beta j)^{-1/2}/2}
        \le C\ex^{-\frac{t^2}{8\beta j}}\, ,
\end{aligned}
\end{equation}}

\noindent
where we have used Lemma \ref{LemmaCompareBridgeFree} and a standard estimate for the maximum of the absolute value of the Brownian motion, see \cite[Equation 4.1.1.4]{borodin1996handbook}.

Choose $t=8j\sqrt{\beta \abs{\mu}}$. Then
\begin{equation}
    \E_{0,0}^{\beta j}\ek{\ex^{-\beta\H(\omega)}\diam(\omega)^\alpha\1\{\diam(\omega)^\alpha>t^\alpha\}}=\Ocal\rk{j^{(1+\alpha)}\ex^{-j\abs{\mu}}}\, ,
\end{equation}
as
\begin{equation}
    \begin{split}
      \E_{0,0}^{\beta j}&\ek{\ex^{-\beta\H(\omega)}\diam(\omega)^\alpha\1\{\diam(\omega)^\alpha>t^\alpha\}}=\int_t^\infty x^\alpha\d\P_{0,0}^{\beta j}(\diam(\omega)\in\d x)\\
      &\le C\int_t^\infty x^\alpha\ex^{-\frac{x^2}{8\beta j}}\d x\le C j^{1/2+\alpha/2}\int_{\sqrt{8j\abs{\mu} }}^\infty x^\alpha \ex^{-x^2}\d x =\Ocal\rk{j^{(1+\alpha)}\ex^{-j\abs{\mu}}}\, ,
    \end{split}
\end{equation}
where we have used a standard estimate for the incomplete Gamma function in the last step.

On the other hand,
\begin{equation}
    \sum_{j\ge 1}\frac{\ex^{\beta \mu j}}{j} \E_{0,0}^{\beta j}\ek{\ex^{-\beta \H(\omega)}\diam(\omega)^\alpha\1\{\diam(\omega)^\alpha\le t^\alpha\}}\le C\sum_{j\ge 1}\frac{\ex^{\mu\beta j}}{j^{1-\alpha}} \E_{0,0}^{\beta j}\ek{\ex^{-\beta \H(\omega)}}\, ,
\end{equation} 
which is finite by Lemma \ref{LemmaSingleLoopDecay}, as $\beta\mu<c_\Phi$.

Hence, we can estimate
\begin{equation}\label{EquationsumfiniteDelta}
     \sum_{j\ge 1}\frac{\ex^{\beta \mu j}}{j} \E_{0,0}^{\beta j}\ek{\ex^{-\beta \H(\omega)}\diam(\omega)^\alpha}\le C \sum_{j\ge 1}(\beta j)^{1+\alpha}\ek{\ex^{-\beta\abs{\mu}}+\ex^{ j(\beta\mu-c_\Phi)}}<\infty\, ,
\end{equation}
Recall the Definition \ref{def:P^H_L}, as $\P_\D^\H\rk{\mathrm{S}_\D>k}=o(1)$ and $\abs{\D}$ finite, by Poisson point process property we have
\begin{equation}\label{eq:finite_moments_S}
    \P_\D^\H\rk{\mathrm{S}_\D>k}=(1+o(1))M^\H_\D[\{\omega:\diam(\omega)>k\}]\le k^{-\alpha}\mathrm{M}_\D^\H\ek{\diam(\omega)^\alpha}\, .
\end{equation}
Note that by translation invariance
\begin{equation}
    \gfrak_n\ek{\1\gk{\mathrm{S}_\D>k}}=\frac{1}{\abs{\L_n}}\sum_{x\in\L_n\cap\Z^d}\gfrak_n\ek{\1\gk{\mathrm{S}_{\D+x}>k}}\, ,
\end{equation}
and hence, by Lemma \ref{lem:FKG},
\begin{equation}
    \gfrak_n\ek{\1\gk{\mathrm{S}_\D>k}}=\delta_{\L_n}^\dir\rk{\1\gk{\mathrm{S}_\D>k}|\emptyset}\le \P_\D^\H\rk{\mathrm{S}_\D>k}\, .
\end{equation}
By Corollary \ref{CorollaryCompact}
\begin{equation}
   \lim_{n\to\infty} \gfrak_n\ek{\1\gk{\mathrm{S}_\D>k}}=\gfrak\ek{\1\gk{\mathrm{S}_\D>k}}\, .
\end{equation}
Hence
\begin{equation}
    \gfrak\ek{\mathrm{S}_\D^{\alpha-1}}\le C\sum_{k\ge 0}k^{\alpha-2}\gfrak[\1\{\mathrm{S}_\D>k\}]
    \le
    C\sum_{k\ge 0}
    k^{-2}
    <\infty\, .
\end{equation}
This concludes the proof as $\alpha$ can be chosen arbitrarily.
\end{proof}
\subsection{Tempered configurations}\label{sec:4.6}
In order to prove that the limiting measure $\gfrak$ satisfies the DLR equations, we need to show that too irregular configurations $\eta$ have zero mass. 
We need to introduce the concept of \textit{tempered configurations}.
\begin{definition}
For $\alpha>0$ and $K,L>0$, we introduce
\begin{equation}
    \Mcal_{K,L}^\alpha=\gk{\eta\in\Omega\colon \mathrm{N}_{\L_n}(\eta)\le K\abs{\L_n}\textnormal{ and }\sup_{\omega\in\eta_{\L_n}}\diam(\omega)\le n^\alpha+L,\,\forall n\ge 1}.
\end{equation}

We furthermore set $\Mcal^\alpha=\cup_{K,L}\Mcal_{K,L}^\alpha$.
\end{definition}

We remark that in particular, $\sup_{\omega\in\eta_{\L_n}}\vol(\omega)\le {\abs{\L_n}^\alpha}+L$ implies that loops cannot extend to the distance ${\abs{\L_n}^\alpha}+L+1$ away from $\L_n$.

\begin{lemma}\label{LemmaGfrakNiceConfiguration}
For any $\alpha>0$, we have that
\begin{equation}
    \gfrak\rk{\Mcal^\alpha}=1\, .
\end{equation}
\end{lemma}
\begin{proof}
The proof follows the same methods as employed in \cite[Lemma 3.2]{dereudre2009existence}. By Lemma \ref{LemmaConvergenceParticlePower}, $\gfrak\ek{\mathrm{N}_{U}}$ is clearly finite, thus by \cite[Theorem 3.7]{nguyen1979ergodic}, $|\L_n|^{-1}\mathrm{N}_{\L_n}$ is $\gfrak$-almost surely finite. In other words, for $\gfrak$-almost surely every $\eta$, we can find $K$ such that the first condition of $\Mcal_{K,L}^\alpha$ is satisfied. 

For the claim regarding the growth of $\vol$, let $U=[0,1]^d$ and notice that by translation invariance and Lemma \ref{lemmacontroldiam}, for every $\alpha>0$,
\begin{equation}
    \begin{split}
        \sum_{x\in\Z^d}\gfrak\left[\sup_{\omega\in\eta_{x+U}}S(\omega)>\abs{x}^\alpha\right]\le 
        \sum_{x\in\Z^d}\gfrak\left[\mathrm{S}_U>\abs{x}^\alpha\right]
        \le C\gfrak\ek{\mathrm{S}_U^{\frac{d-1}{\alpha}}}<\infty\, .
    \end{split}
\end{equation}
Thus, by the Borel--Cantelli lemma, $\gfrak$-almost surely, 
there exists $L>0$ such that for all $x\in\mathbb Z^d$,
\begin{align}
\sup_{\omega\in\eta_{x+U}}S(\omega)\le\abs{x}^\alpha+L,
\end{align}
and in particular,
\begin{align}
\sup_{\omega\in\eta_{\L_n}}S(\omega)=\sup_{x\in\L_n\cap\mathbb Z^d}\sup_{\omega\in\eta_{x+U}}S(\omega)\le \abs{n}^\alpha+L.
\end{align}

\end{proof}

The next lemma gives a uniform control of the $\gfrak_n$.
\begin{lemma}\label{lem:Mcal}
For every $\alpha\in(0,1]$ and every $\e>0$, there exists  $K,L>0$ such that for all $n\in\N$
\begin{equation}
    \gfrak_n\rk{\Mcal^\alpha_{K,L}}\ge 1-\e\, .
\end{equation}
\end{lemma}
\begin{proof}
By Lemma \ref{lem:FKG} and the FKG inequality,
\begin{equation}
    \gfrak_n\rk{\Mcal^\alpha_{K,L}}\ge \P^\H_{\L_n}\rk{\Mcal^\alpha_{K,L}}\ge\P^\H_{\L_n}\rk{\forall m,\mathrm{N}_{\L_m}\le K\abs{\L_m}}\P^\H_{\L_n}\rk{\forall m,\mathrm{S}_{\L_m}\le m^\alpha+L}\, .
\end{equation}
It then suffices to show that uniformly in $n$,
\begin{align}
&\P^\H_{\L_n}\rk{\exists m,\mathrm{N}_{\L_m}> K\abs{\L_m}}\rightarrow 0,\,K\rightarrow\infty\,,\label{eq:Nvanish}\\
&\P^\H_{\L_n}\rk{\exists m,\mathrm{S}_{\L_m}> m^\alpha+L}\rightarrow 0,\,L\rightarrow\infty\,.\label{eq:Svanish}
\end{align}
For \eqref{eq:Nvanish}, by Markov's inequality,
\begin{align}
\P^\H_{\L_n}\rk{\exists m,\mathrm{N}_{\L_m}> K\abs{\L_m}}
\le\sum_{m=1}^n\P^\H_{\L_n}\rk{\ex^{c_\Phi \mathrm{N}_{\L_m}}> \ex^{c_\Phi K\abs{\L_m}}}
\le\sum_{m=1}^n\P^\H_{\L_n}\ek{\ex^{c_\Phi \mathrm{N}_{\L_m}}}{\ex^{-c_\Phi K\abs{\L_m}}}\,.
\end{align}
Then by the same calculation as in \eqref{eq:finite_exponential_N}, there exists a constant $C>0$ so that 
\begin{align}
\P^\H_{\L_n}\ek{\ex^{c_\Phi \mathrm{N}_{\L_m}}}{\ex^{-c_\Phi K\abs{\L_m}}}
=\rk{\mathbb E^\H_{U}\ek{\ex^{c_\Phi \mathrm{N}_{U}}}}^{|\L_m|}{\ex^{-c_\Phi K\abs{\L_m}}}
\le \ex^{(C-c_\Phi K)|\L_m|}\,.
\end{align}
Then \eqref{eq:Nvanish} follows by summing over $m$ and taking $K\rightarrow\infty$.

The second equation \eqref{eq:Svanish} is similar. Indeed, for every $\theta>0$,
\begin{equation}
\begin{aligned}
\P^\H_{\L_n}\rk{\exists m,\mathrm{S}_{\L_m}> m^\alpha+L}
&\le\sum_{m=1}^n\P^\H_{\L_n}\rk{\mathrm{S}_{\L_m}> m^\alpha+L}\\
&\le\sum_{m=1}^n|\L_m|\P_{U}^\H\rk{\mathrm{S}_{U}> m^\alpha+L}
\le\sum_{m=1}^n m^d\frac{\mathbb E^\H_{U}[\mathrm{S}_U^\theta]}{\rk{m^\alpha+L}^\theta}\,.
\end{aligned}
\end{equation}
We fix some $\theta>d/\alpha$ and use \eqref{eq:finite_moments_S} to deduce that the above sum is uniformly bounded for all $n$. Then by dominated convergence, as $L\rightarrow\infty$ it converges to $0$, concluding the proof of \eqref{eq:Svanish}.
\end{proof}
\begin{remark}
Notice that Lemma \ref{lem:Mcal} is slightly stronger than \cite[(3.16)]{dereudre2009existence}. Heuristically, this is possible because our Hamilton has the additive property
\begin{equation}
\H_\L(\eta)\ge\sum_{\omega\in\eta_\L}\H_\L(\omega),\,
\end{equation}
which allows us to define the intermediate measure $\P^\H_\L$ in Lemma \ref{lem:FKG} and take advantage of the FKG inequality.
\end{remark}

We now show that the partition functions are almost surely finite.
\begin{lemma}
    For every $\D\subset\R^d$ compact
    \begin{equation}
        \gfrak\ek{\mathrm Z_\D^\dir\rk{\eta_\D^c}=0}=\gfrak\ek{\Zf_\D\rk{\eta_\D^c}=0}=\gfrak\ek{\mathrm Z^\mathrm{exc}_\D\rk{\eta_\D^c}=0}=0\, .
    \end{equation}
\end{lemma}
\begin{proof}
    Since $\Zf_\D\rk{\eta_\D^c}\ge \P_\D\rk{\xi=\emptyset}$, the claim for the free model follows immediately. 

    Note that
    \begin{equation}
        \mathrm Z_\D^\dir\rk{\eta_\L^c}\ge \ex^{-\H_\D\rk{\eta_\D^{\dir,c}}}\P_\D^\dir\rk{\xi_\D=\emptyset}\, .
    \end{equation}
    Since $\ex^{-\H_\D\rk{\eta_\D^{\dir,c}}}$ is a local function, we have that
    \begin{equation}
        \gfrak\rk{\H_\D\rk{\eta_\L^{\dir,c}}=\infty}=\lim_{n\to\infty}\gfrak_n\rk{\H_\D\rk{\eta_\L^{\dir,c}}=\infty}\, .
    \end{equation}   
    However, the right hand side is zero by Lemma \ref{Lem:PartFiniteFinite}.
    
    Assume that there is $\delta>0$ such that $\gfrak\rk{\mathrm Z_\D^\mathrm{exc}=0}>\delta$. Then, by Lemma \ref{LemmaGfrakNiceConfiguration}, we can find $K,L>0$ such that $\gfrak\rk{\mathrm Z_\D^\mathrm{exc}=0,\Mcal_{K,L}^\alpha}>\delta/2$. However, that means for some $M=M(K,L,\D)$, we have that
    \begin{equation}
        \gfrak\rk{\mathrm Z_\D^\mathrm{exc}\rk{\eta_{\L_M}}=0}>\delta/4\, ,
    \end{equation}
    with positive probability. As this is a local event, there exists $N>M$ such that for all $n\ge N$,
    \begin{equation}
         \gfrak_n\rk{\mathrm Z_\D^\mathrm{exc}\rk{\eta_{\L_M}}=0}>\delta/8\, .
    \end{equation}
    However, $\mathrm Z_\D^\mathrm{exc}\rk{\eta_{\L_M}}\ge \mathrm Z_\D^\mathrm{exc}\rk{\eta_{\L_n}}$ and hence we arrive at contradiction by Lemma \ref{Lem:ExcPartFiniteFinite}.
\end{proof}
\subsection{DLR-Equations}\label{sec:4.7}
Given the results from the previous sections, proving that the DLR equations hold is fairly straight forward. We follow the standard method laid out in \cite{georgii1988gibbs} and adopted for marked processes in \cite{preston2006random,dereudre2009existence}. For any local and bounded function $f$, we write 
\begin{equation}\label{EquationGibbsProperpty0}
    f_\L(\eta)=\int f(\xi) \delta^\dir_\L(\d \xi|\eta)\, .
\end{equation}
Our goal is to prove:
\begin{theorem}\label{theoremGibbsProperty}
For $f$ as above and $\L\subseteq\R^d$ compact
\begin{equation}\label{EquationGibbsProperpty}
    \gfrak\ek{f}=\gfrak\ek{f_\L}\, .
\end{equation}
\end{theorem}
{
Before proving the above theorem, we show how we can deduce Theorem \ref{THM-Main} from it. Define
\begin{equation}\label{Equationh1}
    \gfrak\hk{1}=\int_{U}\gfrak\circ\tau_x\d x\, ,
\end{equation}
where $U=[0,1]^d$. As $\gfrak$ is translation invariant with respect to translations on $\Z^d$, $\gfrak\hk{1}$ is $\R^d$ translation invariant. Furthermore, note that for any function $f$ that is bounded and local, by the translation invariance of the Gibbs kernel
\begin{equation}
\begin{split}
\gfrak^{(1)}[f_\L]
\stackrel{\eqref{Equationh1}}=\int_U \gfrak\circ\tau_x[f_\L]\d x
=\int_U \gfrak[f_\L\circ\tau_{-x}]\d x
&\stackrel{{\eqref{EquationGibbsProperpty0}}}=\int_U \d x \int \gfrak(\d \eta) \int \delta^\dir_\L(\d \xi|\tau_{-x}(\eta))f(\xi)\\
&=\int_U \d x \int \gfrak(\d \eta) \int \delta^\dir_{\tau_x(\L)}(\d \xi|\eta)f(\tau_{-x}(\xi))\\
&\stackrel{{\eqref{EquationGibbsProperpty0}}}=\int_U \gfrak[(f\circ\tau_{-x})_{\tau_x(\L)}]\d x\\
&\stackrel{{\eqref{EquationGibbsProperpty}}}=\int_U \gfrak[f\circ\tau_{-x}]\d x=\gfrak^{(1)}[f]
\end{split}
\end{equation}
and hence $\gfrak\hk{1}$ is Gibbs. 
}

To prove Theorem \ref{theoremGibbsProperty} we need some preparatory lemmas. Fix $\L$ bounded for the rest of the section. We define $\bar\gfrak_n$ another approximation to $\gfrak$ which, in contrary to $\gfrak_n$, satisfies the Gibbs property. The price we pay for that is that $\bar\gfrak_n$ is no longer a probability measure:
\begin{equation}\label{eq:bar_gfrak}
    \bar\gfrak_n=\frac{1}{n^d}\sum_{\mycom{x\in\L_n\cap\Z^d}{\L\subset\tau_x(\L_n)}}\hat{\gfrak}_n\circ\tau_{-x}\, .
\end{equation}
\begin{lemma}\label{LemmaBarGfrakDLR}
The measure $\bar\gfrak_n$ is a sub-probability measure. Furthermore, for any bounded, local, measurable function $f$ and any $\L\subset\L_n$, 
\begin{equation}
    \bar\gfrak_n\ek{f}=\bar\gfrak_n\ek{f_\L}\, .
\end{equation}
\end{lemma}
\begin{proof}
The first claim is immediate. For the second claim, by \eqref{EquationGibbsProperpty0} and \eqref{eq:bar_gfrak},
\begin{equation}
    \begin{split}
        \bar\gfrak_n\ek{f_\L}
        &=\frac{1}{n^d}\sum_{\mycom{x\in\L_n\cap\Z^d}{\L\subset\tau_x(\L_n)}}\int f_\L(\eta)\hat{\gfrak}_n\circ\tau_{-x}(\d \eta)\\
        &=\frac{1}{n^d}\sum_{\mycom{x\in\L_n\cap\Z^d}{\L\subset\tau_x(\L_n)}}\int f_\L(\tau_x(\eta))\hat{\gfrak}_n(\d \eta)\\
        &=\frac{1}{n^d}\sum_{\mycom{x\in\L_n\cap\Z^d}{\L\subset\tau_x(\L_n)}}\int f(\xi)\delta_\L^\dir(\d\xi|\tau_x(\eta))\hat{\gfrak}_n(\d \eta)\\
        &=\frac{1}{n^d}\sum_{\mycom{x\in\L_n\cap\Z^d}{\L\subset\tau_x(\L_n)}}\int f(\tau_x(\xi))\delta^\dir_{\tau_{-x}(\L)}(\d \xi|\eta)\hat\gfrak_n(\d \eta)\,,
    \end{split}
\end{equation}
where the last step is illustrated in Figure \ref{fig:tau}.
\begin{figure}[h]
    \centering
    \includegraphics[width=0.55\textwidth]{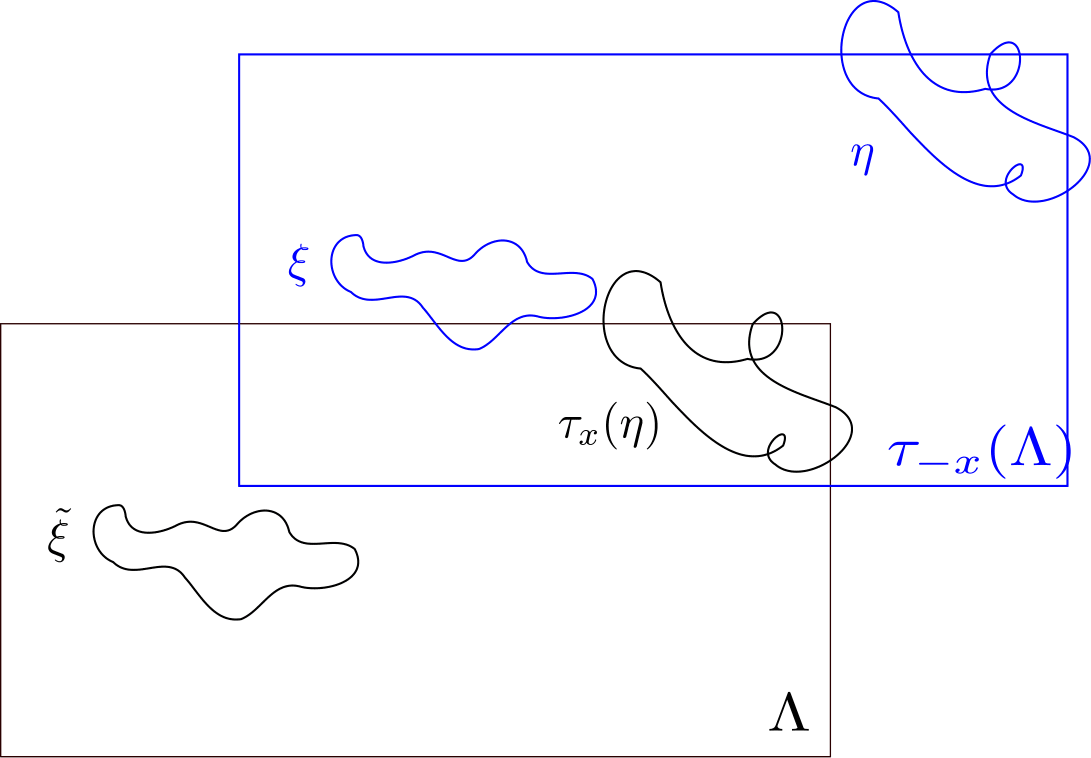}
    \caption{An illustration of $f(\tilde\xi)\delta_\L^\dir(\d\tilde\xi|\tau_x(\eta))=f(\tau_x(\xi))\delta_{\tau_{-x}(\L)}^\dir(\d\xi|\eta)$. The left hand side is shown in black and right hand side in blue.
    }
    \label{fig:tau}
\end{figure}
Since $\L\subset\tau_x(\L_n)$, by Lemma \ref{LemmaConsistent},
we have $\delta^\dir_{\tau_{-x}(\L)}\hat\gfrak_n=\hat\gfrak_n$.
Therefore,
\begin{equation}
    \begin{split}
        \bar\gfrak_n\ek{f_\L}&=\frac{1}{n^d}\sum_{\mycom{x\in\L_n\cap\Z^d}{\L\subset\tau_x(\L_n)}}\int f(\tau_x(\xi))\delta^\dir_{\tau_{-x}(\L)}(\d \xi|\eta)\hat\gfrak_n(\d \eta)\\
        &=\frac{1}{n^d}\sum_{\mycom{x\in\L_n\cap\Z^d}{\L\subset\tau_x(\L_n)}}\int f(\tau_x(\xi))\hat\gfrak_n(\d \xi)\\
        &=\frac{1}{n^d}\sum_{\mycom{x\in\L_n\cap\Z^d}{\L\subset\tau_x(\L_n)}}\int f(\xi)\hat\gfrak_n\circ\tau_{-x}(\d \xi)=\bar\gfrak_n[f]\, .
    \end{split}
\end{equation}

\end{proof}
Next, we show the asymptotic equivalence between $\gfrak$ and $\bar\gfrak_n$.
\begin{lemma}\label{LemmaAsympEquiv}
Let $\psi$ be a function satisfying the condition from Proposition \ref{PropTopology}, then for any $F$ which is local and $\psi$-tame
\begin{equation}
    \lim_{n\to\infty}\Babs{\gfrak_n\ek{F}-\bar\gfrak_n\ek{F}}=0\, .
\end{equation}
\end{lemma}
\begin{proof}
Set $\D\subseteq\R^d$ the support of $F$ and assume that $n$ is sufficiently large for $\Delta\subset\L_n$ and also $\L\subset\L_n$. We then have, using the definition of $\psi$-tame,
\begin{equation}
    \Babs{\gfrak_n\ek{F}-\bar\gfrak_n\ek{F}}\le\frac{1}{n^d}\sum_{\mycom{x\in\L_n\cap\Z^d}{\D\cup\L\nsubseteq\tau_x(\L_n)}}\Babs{\int F(\eta)\,\hat\gfrak_n\circ\tau_{-x}(\d \eta)}\le\frac{C}{n^d}\sum_{\mycom{x\in\L_n\cap\Z^d}{\D\cup\L\nsubseteq\tau_x(\L_n)}}\Babs{\int \left(1+\mathrm{N}_\D^\psi(\eta)\right)\hat\gfrak_n\circ\tau_{-x}(\d \eta)}\, .
    \end{equation}
Note that there exists a constant $C_{\D,\L}$ such that
\begin{equation}
    \#\gk{x\in\L_n\cap\Z^d\colon \Delta\cup\L\nsubseteq\tau_x(\L_n)}\le C_{\D,\L}n^{d-1}\, .
\end{equation}
Thus,
\begin{equation}
    \Babs{\gfrak_n\ek{F}-\bar\gfrak_n\ek{F}}\le C\frac{C_{\D,\L}n^{d-1}}{n^d}+C\frac{C_{\D,\L}n^{d-1}}{n^d}\sup_{x\in\mathbb Z^d}{\int \mathrm{N}_{\tau_x(\D)}^\psi(\eta)\hat\gfrak_n(\d \eta)}\, .
\end{equation}
By Remark \ref{Rmk:Equationbtoinfty_hat}, we find that $\sup_{n\in\mathbb N}\sup_{x\in\mathbb Z^d}\hat\gfrak_n\ek{\mathrm{N}_{\tau_x(\D)}^\psi}<\infty$, and thus the claim follows.
\end{proof}
We now introduce truncated versions of $Z_\L^\dir(\eta)$ and $f_\L(\eta)$, by only taking into accounts loop started in a bounded region $\L'$ containing $\L$: for any $\L'\supset\L$ bounded, we define
\begin{equation}\label{eq:ZLL}
    \mathrm{Z}_{\L,\L'}^\dir\rk{\eta}=\int \ex^{-\beta\H_\L\rk{\eta_{\L'}-\eta_{\L}^\dir+\xi}+\beta\mu \mathrm{N}_\L\rk{\eta_{\L}-\eta_\L^\dir+\xi}}\d\P_\L^\dir\rk{\xi}\, ,
\end{equation}
and
\begin{equation}\label{eq:gL'}
    f_{\L,\L'}(\eta)=\frac{1}{\mathrm{Z}_{\L,\L'}^\dir\rk{\eta}}\int f\rk{\eta_{\L'}-\eta_{\L}^\dir+\xi}\ex^{-\beta\H_\L\rk{\eta_{\L'}-\eta_{\L}^\dir+\xi}+\beta\mu \mathrm{N}_\L\rk{\eta_{\L}-\eta_\L^\dir+\xi}}\d\P_\L^\dir\rk{\xi}\, .
\end{equation}
By definition, $f_{\L,\L'}(\eta)$ is local and bounded.
In the next lemma, we prove that for sufficiently nice $\eta$, the above truncations are sufficiently accurate.

\begin{lemma}\label{LemmaTtruncations}
For any $\alpha \in(0,1)$ and any $\e,K,L>0$, we can find $\L'$ sufficiently large, such that
\begin{equation}
    \sup_{\eta\in\Mcal_{K,L}^\alpha}\babs{f_{\L,\L'}(\eta)-f_\L(\eta)}<\e\, .
\end{equation}
\end{lemma}
\begin{proof}
We first show that for all $\eta\in\Mcal_{K,L}^\alpha$,
\begin{equation}\label{eq:ZZgoal1}
    \babs{\mathrm{Z}_{\L,\L'}^\dir(\eta)-\mathrm{Z}^\dir_\L(\eta)}<\e\, .
\end{equation}
For any $m<n$, take two loops $\omega'\in\eta_{\L_n\setminus\L_{n-1}}$ and $\omega\in\eta_{\L_m}$, then because $\ell(\omega)\le m^\alpha+L$ and $\ell(\omega')\le n^\alpha+L$ (see Figure \ref{fig:L_m_and_L_n}),
\begin{equation}\label{EquationHelpEstimateDifference}
    \inf_{\mycom{0\le t\le \beta \ell(\omega')}{0\le s\le \beta \ell(\omega)}}\abs{\omega'(t)-\omega(s)}\ge n-1-n^\alpha-L-1-m-m^\alpha-L-1\, .
\end{equation}
Therefore, if $m<n$ are large enough, we have by Assumption \ref{refassump}
{\begin{equation}\label{EquationHelpEstimateDifference2}
T(\omega,\omega')\le \ell(\omega)\ell(\omega')\Psi(n-n^\alpha-m-m^\alpha-2L-3)\, .
\end{equation}}

\begin{figure}[h]
    \centering
    \includegraphics[width=0.55\textwidth]{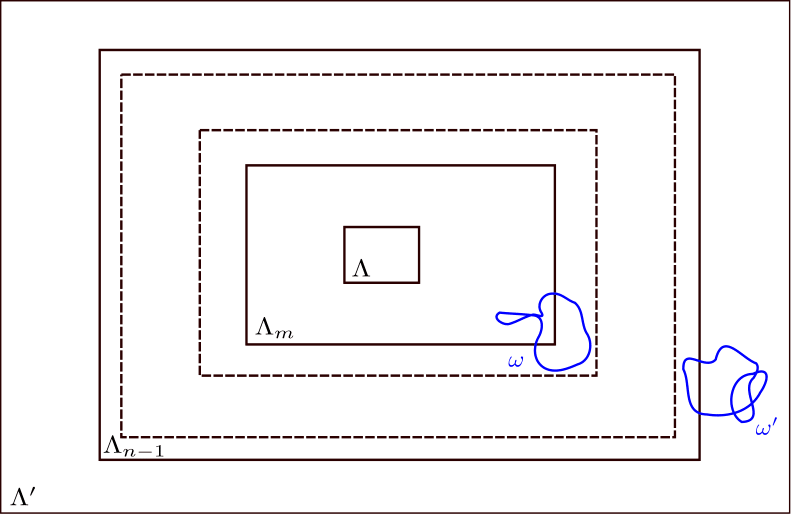}
    \caption{An illustration of $\ell(\omega)\le m^\alpha+L$ and $\ell(\omega')\le n^\alpha+L$. As the radius of the two loops are bounded (not reaching the dotted lines), we have \eqref{EquationHelpEstimateDifference}.
    }
    \label{fig:L_m_and_L_n}
\end{figure}

By decomposing $\eta=\eta_{\L'}+\eta_{\L'}^c$, we get that
\begin{equation}\label{EquationDifferenceExpanded220}
    \babs{\mathrm{Z}_{\L,\L'}^\dir(\eta)-\mathrm{Z}^\dir_\L(\eta)}=\int\ex^{-\beta\H_\L\rk{\eta_{\L'}-\eta_\L^\dir+\xi}+\beta\mu \mathrm{N}_\L(\eta_{\L}-\eta_\L^\dir+\xi)}\Babs{\ex^{-\sum_{\omega \in \xi+\eta_{\L}-\eta_{\L}^\dir}\sum_{\omega'\in\eta_{\L'}^c}T(\omega,\omega')}-1}\d\P_\L^\dir(\xi)\, .
\end{equation}
Without loss of generality, we can choose large enough $\L'$, $m$ and $n-m$ such that $\L\subset\L_m\subset\L_n\subset\L'$. We bound
\begin{equation}\label{eq:ZZbound}
        \eqref{EquationDifferenceExpanded220}\le\int\ex^{-\beta\H_\L(\xi)+\max\{\mu K|\L|,0\}+\beta\mu \mathrm{N}_\L(\xi)}\left(1-\ex^{-\sum_{\omega \in \xi+\eta_{\L}-\eta_{\L}^\dir}\sum_{\omega'\in\eta_{\L'}^c}T(\omega,\omega')}\right)\d\P_\L^\dir(\xi)\, ,
\end{equation}
where we used the fact that $\H_\L\rk{\eta_{\L'}-\eta_\L^\dir+\xi}\ge \H_\L(\xi)$ and $\mathrm{N}_\L(\eta_{\L}-\eta_\L^\dir+\xi)=\mathrm{N}_\L(\eta_{\L}-\eta_\L^\dir)+\mathrm{N}_\L(\xi)\le K\abs{\L}+\mathrm{N}_\L(\xi)$.
Abbreviate $C_o=m+m^\alpha+2L+3$. Write $o_n(1)$ for a null sequence as $n\to\infty$ and note that by using Equation \eqref{EquationHelpEstimateDifference2}
\begin{equation}\label{eq:ZZbound2}
    \begin{split}
        \sum_{\omega \in \xi+\eta_{\L}-\eta_{\L}^\dir}\sum_{\omega'\in\eta_{\L'}^c}T(\omega,\omega')&\le \mathrm{N}_\L(\xi+\eta_{\L}-\eta_{\L}^\dir)\sum_{k> n}\mathrm{N}_{\L_k\setminus \L_{k-1}}(\eta_{\L'}^c) |\Psi\rk{k-k^\alpha-C_o}|\\
        &\le C K \mathrm{N}_\L(\xi+\eta_{\L}-\eta_{\L}^\dir)\sum_{k> n}k^{d-1}|\Psi\rk{k-k^\alpha-C_o}|\\
        &\le C K \mathrm{N}_\L(\xi+\eta_{\L}-\eta_{\L}^\dir) o_n(1)\\
        &\le C K  o_n(1)(\mathrm{N}_\L(\xi)+K|\L|)\, ,
    \end{split}
\end{equation}
where we used the volume growth of $\mathrm{N}_{\L_n}(\eta)$ under the assumption that $\eta\in\Mcal_{K,L}^\alpha$ and the integrability of $x^{d-1}\Psi(x)$.
Therefore
\begin{equation}
    \begin{split}
        \babs{\mathrm{Z}_{\L,\L'}^\dir(\eta)-\mathrm{Z}^\dir_\L(\eta)}&\le\int\ex^{-\beta\H_\L(\xi)+\max\{\mu K|\L|,0\}+\beta\mu \mathrm{N}_\L(\xi)}\left(1-\ex^{-CKo_n(1)(\mathrm{N}_\L(\xi)+K|\L|)}\right)\d\P_\L^\dir(\xi)\,,
    \end{split}
\end{equation}
which is finite for some $n$ large enough (see Equation \eqref{Integrabilityofhamiltonian}), and converges to $0$ as $n\rightarrow\infty$ by dominated convergence.
Thus \eqref{eq:ZZgoal1} is proved.

For the approximation of $f_\L$ by $f_{\L,\L'}$, we can use a similar procedure. To increase readability, we abbreviate
\begin{equation}
    \widetilde{\H}(\eta,\xi)=\beta\H_\L\rk{\eta-\eta_{\L}^\dir+\xi}\quad\text{and}\quad\widetilde{\H}_{\mathrm{trunc}}(\eta,\xi)=\beta\H_\L\rk{\eta_{\L'}-\eta_{\L}^\dir+\xi}\, ,
\end{equation}
as well as
\begin{equation}
    \widetilde{\mathrm{N}}(\eta,\xi)=\beta\mu \mathrm{N}_\L(\eta_\L-\eta^\dir_\L+\xi)\, .
\end{equation}
Indeed, 
we expand
\begin{equation}
\begin{aligned}
&|f_\L-f_{\L,\L'}|\\
&=\left|
\frac{1}{\mathrm{Z}_{\L}^\dir\rk{\eta}}\int f\rk{\eta-\eta_{\L}^\dir+\xi}\ex^{-\widetilde{\H}(\eta,\xi)+\widetilde{\mathrm{N}}(\eta,\xi)}\d\P_\L^\dir\rk{\xi}\right.\\
&\qquad\qquad\qquad\qquad\qquad\qquad\qquad-\left.
\frac{1}{\mathrm{Z}_{\L,\L'}^\dir\rk{\eta}}\int f\rk{\eta_{\L'}-\eta_{\L}^\dir+\xi}\ex^{-\widetilde{\H}_{\mathrm{trunc}}(\eta,\xi)+\widetilde{\mathrm{N}}(\eta,\xi)}\d\P_\L^\dir\rk{\xi}
\right|\\
&\le\frac{1}
{\mathrm{Z}_{\L}^\dir\rk{\eta}}\int \Babs{
f\rk{\eta-\eta_{\L}^\dir+\xi}\ex^{-\widetilde{\H}(\eta,\xi)}- f\rk{\eta_{\L'}-\eta_{\L}^\dir+\xi}\ex^{-\widetilde{\H}_{\mathrm{trunc}}(\eta,\xi)}}\ex^{\widetilde{\mathrm{N}}(\eta,\xi)}\d\P_\L^\dir\rk{\xi}\\
&\qquad\qquad\qquad\qquad\qquad\qquad\qquad +\frac{\Babs{\mathrm{Z}_{\L}^\dir\rk{\eta}-\mathrm{Z}^\dir_{\L,\L'}\rk{\eta}}}{\mathrm{Z}^\dir_{\L}\rk{\eta}\mathrm{Z}^\dir_{\L,\L'}\rk{\eta}}{\int f\rk{\eta_{\L'}-\eta_{\L}^\dir+\xi}\ex^{-\widetilde{\H}_{\mathrm{trunc}}(\eta,\xi)+\widetilde{\mathrm{N}}(\eta,\xi)}\d\P_\L^\dir\rk{\xi}}\,.
\end{aligned}
\end{equation}
Let $\sup_\eta |f(\eta)|=C_f<\infty$, then $|f_\L-f_{\L,\L'}|$ can be further bounded above by
\begin{equation}
\begin{aligned}
\frac{C_f}
{\mathrm{Z}_{\L}^\dir\rk{\eta}}&\int \Babs{
\ex^{-\widetilde{\H}(\eta,\xi)}- \ex^{-\widetilde{\H}_{\mathrm{trunc}}(\eta,\xi)}}\ex^{\widetilde{\mathrm{N}}(\eta,\xi)}\d\P_\L^\dir\rk{\xi}
+
\frac{C_f\Babs{\mathrm{Z}_{\L}^\dir\rk{\eta}-\mathrm{Z}^\dir_{\L,\L'}\rk{\eta}}}{\mathrm{Z}_{\L}^\dir\rk{\eta}}\\
&=2C_f\frac{\int 
\rk{
\ex^{-\widetilde{\H}_{\mathrm{trunc}}(\eta,\xi)}
-
\ex^{-\widetilde{\H}(\eta,\xi)}
}
\ex^{\widetilde{\mathrm{N}}(\eta,\xi)}\d\P_\L^\dir\rk{\xi}}
{\int
\ex^{-\widetilde{\H}(\eta,\xi)+\widetilde{\mathrm{N}}(\eta,\xi)}\d\P_\L^\dir\rk{\xi}}\\
&=2C_f\frac{\int 
\ex^{\beta\H_\L(\eta-\eta_\L^\dir)-\widetilde{\H}(\eta,\xi)+\widetilde{\mathrm{N}}(\eta,\xi)}
\rk{
\ex^{\sum_{\omega \in \xi+\eta_{\L}-\eta_{\L}^\dir}\sum_{\omega'\in\eta_{\L'}^c}T(\omega,\omega')}
-1
}
\d\P_\L^\dir\rk{\xi}}
{\int
\ex^{\beta\H_\L(\eta-\eta_\L^\dir)-\widetilde{\H}(\eta,\xi)+\widetilde{\mathrm{N}}(\eta,\xi)}\d\P_\L^\dir\rk{\xi}}\,.
\end{aligned}
\end{equation}
Since we have non-negative Hamilton, we can conclude that
\begin{equation}\label{eq:f_minus_g}
\begin{aligned}
|f_\L-f_{\L,\L'}|&\le
2C_f\frac{\int \ex^{-\beta\H_\L(\xi)+\beta\mu \mathrm{N}_\L(\xi)}
\rk{
\ex^{\sum_{\omega \in \xi+\eta_{\L}-\eta_{\L}^\dir}\sum_{\omega'\in\eta_{\L'}^c}T(\omega,\omega')}
-1
}
\d\P_\L^\dir\rk{\xi}}
{\int
\ex^{\beta\H_\L(\eta-\eta_\L^\dir)-\widetilde{\H}(\eta,\xi)+\beta\mu \mathrm{N}_\L(\xi)}\d\P_\L^\dir\rk{\xi}}\,.
\end{aligned}
\end{equation}

Now the denominator is uniformly bounded below by
\begin{equation}\label{eq:denominator}
\int
\ex^{\beta\H_\L(\eta-\eta_\L^\dir)-\widetilde{\H}(\eta,\xi)+\beta\mu \mathrm{N}_\L(\xi)}\d\P_\L^\dir\rk{\xi}\ge \P^\dir_\L(\emptyset)>0\,,
\end{equation}
and for the numerator we use \eqref{eq:ZZbound2} to get
\begin{equation}
\begin{split}
\int &
\ex^{-\beta\H_\L(\xi)+\beta\mu \mathrm{N}_\L(\xi)}\rk{
\ex^{\sum_{\omega \in \xi+\eta_{\L}-\eta_{\L}^\dir}\sum_{\omega'\in\eta_{\L'}^c}T(\omega,\omega')}-1
}
\d\P_\L^\dir\rk{\xi}\\
&\le\int
\ex^{-\beta\H_\L(\xi)+\beta\mu \mathrm{N}_\L(\xi)}
\rk{\ex^{C K  o_n(1)(\mathrm{N}_\L(\xi)+K|\L|)}-1}\d\P_\L^\dir\rk{\xi},
\end{split}
\end{equation}
which is finite for some $n$ large enough, and converging to $0$ as $n\rightarrow\infty$ by dominated convergence.

In summary, we notice that the bounds above are uniform in $\eta\in\Mcal^\alpha_{K,L}$, and
we can take $n$ and $\L'\supset\L_n$ large enough so that $\sup_{\eta\in\Mcal_{K,L}^\alpha}\Babs{f_\L(\eta)-f_{\L,\L'}(\eta)}$ is arbitrarily small.
\end{proof}
We are now in the position to prove Theorem \ref{theoremGibbsProperty}.
\begin{proof}[{Proof of Theorem \ref{theoremGibbsProperty}}]
Fix $\epsilon>0$. By Lemma \ref{LemmaGfrakNiceConfiguration} and Lemma \ref{lem:Mcal}, we can fix $K,L>0$ such that 
\begin{equation}
\gfrak((\Mcal^\alpha_{K,L})^c)\le\epsilon,\quad\bar\gfrak_n((\Mcal^\alpha_{K,L})^c)\le\gfrak_n((\Mcal^\alpha_{K,L})^c)\le \epsilon\,.
\end{equation}
Fix also $\L'$ such that Lemma \ref{LemmaTtruncations} holds.

Since $f$ and $f_{\L,\L'}$ are local and bounded, by Lemma \ref{LemmaAsympEquiv},
\begin{equation}
\Babs{\bar\gfrak_n[f]-\gfrak_n[f]}\rightarrow 0,\quad
\Babs{\bar\gfrak_n[f_{\L,\L'}]-\gfrak_n[f_{\L,\L'}]}\rightarrow 0\,.
\end{equation}
Further, since $|f|,|f_\L|$ and $|f_{\L,\L'}|$ are all bounded by $\sup |f|$, by Lemma \ref{LemmaTtruncations},
\begin{equation}
\Babs{\bar\gfrak_n[f_\L]-\bar\gfrak_n[f_{\L,\L'}]}
\le \sup |f|\cdot\bar\gfrak_n[(\Mcal_{K,L}^\alpha)^c]+\epsilon
\le (1+\sup|f|)\epsilon\,,
\end{equation}
and similarly
\begin{equation}
\Babs{\gfrak[f_\L]-\gfrak[f_{\L,\L'}]}
\le (1+\sup|f|)\epsilon\,.
\end{equation}
In addition, since $f$ and $f_{\L,\L'}$ are local and bounded,
by Proposition \ref{PropTopology},
\begin{equation}
\gfrak_n[f]\rightarrow\gfrak[f],\quad\gfrak_n[f_{\L,\L'}]\rightarrow\gfrak[f_{\L,\L'}]\,.
\end{equation}
Finally, recall that we have
\begin{equation}
    \bar\gfrak_n\ek{f}=\bar\gfrak_n\ek{f_\L}\, 
\end{equation}
by Lemma \ref{LemmaBarGfrakDLR}, we conclude that
\begin{equation}
\begin{aligned}
\Babs{\gfrak[f]-\gfrak[f_\L]}
\le&\Babs{\gfrak[f]-\gfrak[f_{\L,\L'}]}+(1+\sup|f|)\epsilon\\
\le&\limsup_{n\rightarrow\infty}\Babs{\gfrak_n[f]-\gfrak_n[f_{\L,\L'}]}+(1+\sup|f|)\epsilon\\
=&\limsup_{n\rightarrow\infty}\Babs{\bar\gfrak_n[f]-\bar\gfrak_n[f_{\L,\L'}]}+(1+\sup|f|)\epsilon\\
=&\limsup_{n\rightarrow\infty}\Babs{\bar\gfrak_n[f_\L]-\bar\gfrak_n[f_{\L,\L'}]}+(1+\sup|f|)\epsilon\le 2(1+\sup|f|)\epsilon\,.
\end{aligned}
\end{equation}
The conclusion follows as $\epsilon$ can be arbitrarily small.
\end{proof}
To show that $\gfrak$ is also invariant under the application of $\rk{\gamma_\L}_\L$ and $\rk{\df_\L}_\L$, we need a preparatory lemma. See Figure \ref{fig:lemmakernel} for an illustration.
\begin{lemma}\label{LemmaDandLkernelconsistency}
Fix $\L\subset\D$. If $\eta$ is such that all loops intersecting $\L$ are contained in $\D$, then
\begin{equation}
    \delta_\D^\dir\gamma_\L(\eta)=\delta_\D^\dir(\eta)\quad\textnormal{and}\quad \delta_\D^\dir\df_\L(\eta)=\delta_\D^\dir(\eta)\, .
\end{equation}
\end{lemma}
\begin{proof}
We only prove for $\rk{\delta^\exc_\L}$, as the proof for $\rk{\df_\L}$ is similar.
We begin by expanding
\begin{equation}
    \delta_\D^\dir\gamma_\L(A|\eta)=\frac{1}{\mathrm{Z}_\D^\dir\rk{\eta_\D^{\dir,c}}}\int \gamma_\L(A|\eta_\D^{\dir,c}+\zeta)\ex^{-\beta\H_\D(\eta_\D^{\dir,c}+\zeta)+\beta\mu \mathrm{N}_\D(\eta_\D^{\dir,c}+\zeta)}\d\P_\D^\dir(\zeta)\, .
\end{equation}
By our assumption on $\eta$, we have $\eta^\exc_{\L^c}=\rk{\eta_\D^\dir}^\exc_{\L^c}$ and $\Q_\L\rk{\cdot\big|\cdot+\eta_\D^{\dir,c}}=\Q_\L\rk{\cdot\big|\cdot}$, thus
\begin{multline}
    \gamma_\L(A|\eta_\D^{\dir,c}+\zeta)\\
    =\frac{1}{{\mathrm{Z}}^\exc_\L\rk{\eta_\D^{\dir,c}+\zeta^\exc_{\L^c}}}\int\1_A\rk{\eta_\D^{\dir,c}+\xi+\zeta^\exc_{\L^c}}\ex^{-\beta\H_\L\rk{\eta_\D^{\dir,c}+\xi+\zeta^\exc_{\L^c}}+\beta\mu \mathrm{N}_\L(\eta_\D^{\dir,c}+\xi+\zeta^\exc_{\L^c})}\d \Q_\L\rk{\xi\big|\zeta^\exc_{\L^c}}\, ,
\end{multline}
see Figure \ref{fig:lemmakernel}.
\begin{figure}[h]
    \centering
    \includegraphics[width=0.65\textwidth]{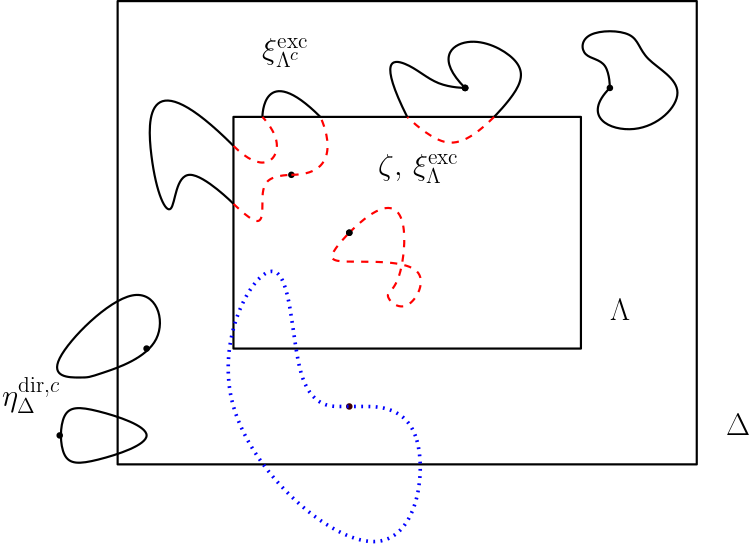}
    \caption{An illustration of the situation in Lemma \ref{LemmaDandLkernelconsistency}: configurations which intersect both $\D$ and $\L^c$ like the dotted, blue one in the above figure are excluded. Hence $\Q_\L$ is a function of the loops started inside $\D$, i.e.,  $\Q_\L\rk{\cdot\big|\eta}=\Q_\L\rk{\cdot\big|\eta_\D^\dir}$}
    \label{fig:lemmakernel}
\end{figure}
By Lemma \ref{lem:H_trick_exc},
\begin{equation}
    \H_\D(\eta_\D^{\dir,c}+\zeta)-\H_\L(\eta_\D^{\dir,c}+\zeta)=\H_\D(\eta_\D^{\dir,c}+\zeta^\exc_{\L^c}+\xi)-\H_\L(\eta_\D^{\dir,c}+\zeta^\exc_{\L^c}+\xi)\, .
\end{equation}
Furthermore, similar to Lemma \ref{lem:consistency_Q}, we
have
\begin{equation}
\d\P^\dir_{\D}(\zeta)
=\d\Q_{\D}(\zeta\big|\varnothing)
=\d\Q_{\L}(\zeta^\exc_\L\big|\zeta)\d\tilde\Q_{\D\setminus\L}(\zeta\big|\varnothing)\,,
\end{equation}
where $\tilde \Q_{\L\setminus\D}(\cdot\big|\varnothing)$ samples loops inside $\D$, and then cuts them at the boundary of $\L$ by taking $(\cdot)^\exc_{\L^c}$.
We deduce analogously to Lemma \ref{LemmaConsistentQkernel}:
{\small\begin{equation}
    \begin{split}
        &\mathrm{Z}_\D^\dir\rk{\eta_\D^{\dir,c}}\delta_\D^\dir\gamma_\L(A|\eta)\\
        =& \int\frac{\1_A\rk{\xi+\zeta^\exc_{\L^c}+\eta_\D^{\dir,c}}}{{\mathrm{Z}}^\exc_\L\rk{\zeta^\exc_{\L^c}+\eta_\D^{\dir,c}}}
        \ex^{-\beta\H_\D\rk{\zeta+\eta_\D^{\dir,c}}-\beta\H_\L(\eta_\D^{\dir,c}+\zeta^\exc_{\L^c}+\xi)
        +\beta\mu \mathrm{N}_\D\rk{\zeta+\eta_\D^{\dir,c}}+\beta\mu \mathrm{N}_\L(\eta_\D^{\dir,c}+\zeta^\exc_{\L^c}+\xi)
        }\d \Q_\L\rk{\xi\big|\zeta}
        \d\Q_\D\rk{\zeta\big|\varnothing}\\
        =&\int\frac{\1_A\rk{\xi+\zeta^\exc_{\L^c}+\eta_\D^{\dir,c}}}{{\mathrm{Z}}^\exc_\L\rk{\zeta^\exc_{\L^c}+\eta_\D^{\dir,c}}}
        \ex^{-\beta\H_\D\rk{\zeta^{\exc}_\L+\zeta^{\exc}_{\L^c}+\eta_\D^{\dir,c}}-\beta\H_\L(\eta_\D^{\dir,c}+\zeta^\exc_{\L^c}+\xi)
        +\beta\mu \mathrm{N}_\D\rk{\zeta^{\exc}_\L+\zeta^{\exc}_{\L^c}+\eta_\D^{\dir,c}}+\beta\mu \mathrm{N}_\L(\eta_\D^{\dir,c}+\zeta^\exc_{\L^c}+\xi)
        }\\
        &\qquad\quad\qquad\qquad\qquad\qquad\qquad\qquad\qquad\qquad\qquad\qquad\qquad\qquad\qquad\qquad
        \d\Q_\L\rk{\zeta^\exc_{\L}\big|\zeta}
        \d \Q_\L\rk{\xi\big|\zeta}
        \d\tilde\Q_{\D\setminus\L}\rk{\zeta\big|\varnothing}\\
        =&\int\frac{\1_A\rk{\xi+\zeta^\exc_{\L^c}+\eta_\D^{\dir,c}}}{{\mathrm{Z}}^\exc_\L\rk{\zeta^\exc_{\L^c}+\eta_\D^{\dir,c}}}
        \Bigg(\ex^{-\beta\H_\L(\eta_\D^{\dir,c}+\zeta^\exc_{\L^c}+\xi)
        +\beta\mu \mathrm{N}_\L(\eta_\D^{\dir,c}+\zeta^\exc_{\L^c}+\xi)
        } \d\Q_\L\rk{\zeta^\exc_{\L}\big|\zeta}\Bigg)\\
        &\qquad\qquad\qquad\qquad\qquad\qquad\qquad\qquad\qquad
        \times\ex^{-\beta\H_\D\rk{\zeta^{\exc}_\L+\zeta^{\exc}_{\L^c}+\eta_\D^{\dir,c}}
        +\beta\mu \mathrm{N}_\D\rk{\zeta^{\exc}_\L+\zeta^{\exc}_{\L^c}+\eta_\D^{\dir,c}}
        }\d \Q_\L\rk{\xi\big|\zeta}
        \d\tilde\Q_{\D\setminus\L}\rk{\zeta\big|\varnothing}\\
        =&\int\frac{\1_A\rk{\xi+\zeta^\exc_{\L^c}+\eta_\D^{\dir,c}}{\mathrm{Z}}^\exc_\L\rk{\zeta^\exc_{\L^c}+\eta_\D^{\dir,c}}}{{\mathrm{Z}}^\exc_\L\rk{\zeta^\exc_{\L^c}+\eta_\D^{\dir,c}}}\ex^{-\beta\H_\D(\eta_\D^{\dir,c}+\zeta^\exc_{\L^c}+\xi)
        +\beta\mu \mathrm{N}_\D(\eta_\D^{\dir,c}+\zeta^\exc_{\L^c}+\xi)}\d \Q_\L\rk{\xi\big|\zeta^\exc_{\L^c}}
        \d\tilde\Q_{\D\setminus\L}\rk{\zeta^\exc_{\L^c}\big|\varnothing}\\
        =&\int{\1_A\rk{\alpha+\eta_\D^{\dir,c}}}\ex^{-\beta\H_\D(\eta_\D^{\dir,c}+\alpha)
        +\beta\mu \mathrm{N}_\D(\eta_\D^{\dir,c}+\alpha)}
        \d \P^\dir_\D(\alpha)\\
        =&\mathrm{Z}_\D^\dir\rk{\eta_\D^{\dir,c}}\delta_\D^\dir(A|\eta)\, .
    \end{split}
\end{equation}}
This concludes the proof.
\end{proof}
We can now prove the DLR-equations in a path-wise sense.
\begin{theorem}\label{TheoremPathDLR}
For $\L$ finite, we have
\begin{equation}
    \gfrak\gamma_\L=\gfrak\quad \textnormal{and}\quad \gfrak\df_\L=\gfrak\, .
\end{equation}
\end{theorem}
\begin{proof}
We only give the proof for the case of $\rk{\gamma_\L}_\L$, as the free kernel is analogous.

By Theorem \ref{theoremGibbsProperty}, we have
\begin{equation}
    \gfrak\gamma_\L=\gfrak\delta_\D\gamma_\L\, ,
\end{equation}
for any $\Delta$. Lemma \ref{LemmaGfrakNiceConfiguration} gives that there exists $K_o,L_o$ such that for every $K>K_o$ and $L>L_o$ we have $\gfrak\rk{\Mcal^\alpha_{K,L}}\ge 1-\e$. Thus, outside of a set of at most $\e$ mass, we can find $\Delta$ such that $\eta_\L\subset\eta_\D^\dir$. However, this means that by Lemma \ref{LemmaDandLkernelconsistency} that $\delta^\dir_\D \delta^\exc_\L(\eta)=\delta_\D^\dir(\eta)$ and thus
\begin{equation}
    \babs{\gfrak-\gfrak\gamma_\L}=\babs{\gfrak-\gfrak\delta_\D^\dir\gamma_\L}\le \e +\babs{\gfrak-\gfrak\delta^\dir_\D}=\e\, ,
\end{equation}
where $\babs{\nu_1-\nu_2}$ is short for the total variational distance $\sup_{A\in\Fcal}\babs{\nu_1(A)-\nu_2(A)}$.
As $\e>0$ was arbitrary, the result follows.
\end{proof}
\subsection{Proof of Proposition \ref{THMProp}}\label{sec:4.8}
The proof of Proposition \ref{THMProp} is now relatively straight forward. 

We first notice that Lemma \ref{lem:FKG} remains valid if we replace $\delta_\L^\dir$ by $\df_\L$.

Then to prove that 
\begin{equation}
    \mathsf{G}\ek{\ex^{\sum_{\omega\in\eta_\D}\psi(\omega)}}<\infty\, ,
\end{equation}
for $\mathsf{G}$ Gibbs with respect to $\rk{\df_\L}_\L$, we have
\begin{equation}
    \mathsf{G}\ek{\ex^{\sum_{\omega\in\eta_\D}\psi(\omega)}}=\mathsf{G}\ek{\df_\D\ek{\ex^{\sum_{\omega\in\eta_\D}\psi(\omega)}|\eta_\D^c}}\le \E_\D^\H\ek{\ex^{\sum_{\omega\in\eta_\D}\psi(\omega)}}\, .
\end{equation}
By the Campbell formula
\begin{equation}
    \E_\D^\H\ek{\ex^{\sum_{\omega\in\eta_\D}\psi(\omega)}}=\exp\rk{\int_\D\d x\sum_{j\ge 1}\frac{\ex^{\beta\mu j}}{j}\E_{x,x}^{\beta j}\ek{\ex^{-\beta \H(\omega)}\rk{\ex^{\psi(\omega)}-1}}}\, .
\end{equation}
Since $\E_{x,x}^{\beta j}\ek{\ex^{-\H(\omega)}}=\Ocal\rk{\ex^{-c_\Phi j}}$ and $\psi(\omega)\le \alpha \ell(\omega)$, the above sum is finite.

The second statement is derived similarly, now \eqref{eq:finite_moments_S} is used to estimate the powers of the diameter with respect to $\P^\H_\D$.

\newpage
\appendix
\section{Frequently used notations}\label{SectionAppendix}

\begin{table}[H]
  \renewcommand{\arraystretch}{1.2}
\begin{tabular}{ |p{1.2cm}||p{5.0cm}|p{6cm}|p{3.4cm}|  }
 \hline
 Symbol & Definition & Explanation &Class\\
 \hline \hline
  $\mu$   &  $\mu\in\R$, usually $\mu\le c_\Phi/\beta$   &Chemical potential&   Model parameter\\
 \hline
 $\beta$   &  $\beta>0$   &Inverse temperature&   Model parameter\\
 \hline
 $\Phi$   &   $\Phi\colon[0,\infty)\to\R\cup\{+\infty\}$  & Interaction potential &  Model parameter \\
\hline
$\eta$    & $\eta=\sum_{\omega}\delta_\omega$    & Loop configuration &  Configuration \\
\hline
$\eta_\L$   &  $\eta=\sum_{\omega}\delta_\omega
 \1\{\omega(0)\in\L\}$   & Loops started in $\L$& Configuration  \\
\hline
$\eta_\L^c$   &  $\eta_\L^c=\eta-\eta_\L$   & Loops started outside $\L$& Configuration  \\
\hline
$\eta_\L^\dir$   &  $\eta=\sum_{\omega}\delta_\omega
 \1\{\omega\subset\L\}$   & Loops contained in $\L$& Configuration  \\
\hline
$\eta_\L^{\dir,c}$   &  $\eta_\L^{\dir,c}=\eta-\eta_\L^\dir$   & Loops not contained in $\L$& Configuration  \\
\hline
$\eta_\L^\exc$    &  See Eq. \eqref{Equationexceta}   & Excursions inside $\L$& Configuration  \\
\hline
$W$   & See Eq. \eqref{EquationDefW}    & Self interaction & Loop function  \\
 \hline
$T$    &  See Eq. \eqref{EquationDefT}   & Pair interaction &  Loop function \\
\hline
$\H_\L$    &  See Eq. \eqref{Equation H}  & Hamiltonian &   Loop function\\
\hline
$U$    &  See Eq. \eqref{Equation U}   & Interaction &   Loop function\\
\hline
$\ell$    &   See above Eq. \eqref{Eq123three} & Loop particle number & Loop function  \\
\hline
$\vol$    &   See Eq. \eqref{definitinvol} & Maximum diameter & Loop function  \\
\hline
$\mathrm{N}_\L$    &   See Eq. \eqref{definitinvol}  & Total particle number &  Loop function  \\
\hline
$\P_x$    &     & Brownian motion, started at $x\in\R^d$&  Measure \\
\hline
$\B_{x,y}^t$    &     & B. bridge from $x$ to $y$ in time $t$& Measure  \\
\hline
$\P_{x,y}^t$   &  $\P_{x,y}^t=p_t(x,y)\B_{x,y}^t$   & Unnormalized bridge measure & Measure  \\
\hline
$\mathrm{M}_\L$    &   $\mathrm{M}_\L=\int\d x\sum_{j\ge 1}\frac{1}{j}\P_{x,x}^{\beta j}$  & Loop measure&  Measure \\
\hline
$\P_\L$    &  See Definition \ref{definitionRefProc}  & Poisson reference process & Measure  \\
\hline
$\P_\L^\dir$    &  See Definition \ref{definitionRefProc}   & Reference process, Dirichlet b.c. & Measure  \\
\hline
$\Q_\L$    &  See Equation \ref{EquationDefinitionOfQ}   & Excursion Kernel & Measure  \\
\hline
$\gfrak_n,\gfrak$    &  See Eq. \eqref{Equation gn}  & Approximation and Gibbs measure & Measure  \\
\hline
$\delta^\dir_\L$   &  See Eq. \eqref{Equation delta}   &(Dirichlet) Gibbs kernel &   Kernel\\
\hline
$\df_\L$   &  See Eq. \eqref{Equation deltafree}  &Free Gibbs kernel &   Kernel\\
\hline
$\gamma_\L$   &  See Eq. \eqref{Equation deltaexc}  &Excursion Gibbs kernel &   Kernel\\
\hline
$\mathrm{Z}^\dir_\L$   &  See Eq. \eqref{Equation Z}   &(Dirichlet) partition function &   Constant\\
\hline
 $\Zf_\L$   &  See Eq. \eqref{Equation Zfree}   &Free partition function &   Constant\\
\hline
$ \mathrm{Z}_\L^\mathrm{exc}$   &  See Eq. \eqref{Equation Zexc}   &Excursion partition function &   Constant\\
\hline
$c_\Phi$   &  See Eq. \ref{LemmaSingleLoopDecay}   & Decay speed of a single loop &   Constant\\
\hline
$I$ & See Eq. \eqref{Equation I} & Specific entropy & Function \\
\hline
\end{tabular}\label{table1} 
\end{table}
\newpage
\bibliography{thoughts}{}
\bibliographystyle{alpha}
\end{document}